\documentclass{conm-p-l} 
\usepackage{amsmath} 
\usepackage{amssymb} 
\usepackage{amscd}

\usepackage{pb-diagram} 
\usepackage{epsf}

\def\oplusinf{\mathop{\oplus}} 
\def\otimesinf{\mathop{\otimes}}
\def\id{{\mathrm{Id}}} 
\def\im{{\mathrm{Im}}} 
\def\dim{{\mathrm{dim}}}
\def\ker{{\mathrm{Ker}}} 
\def\Der{{\mathrm{Der}}}
\def\ham{{\mathrm{Ham}}} 

\def\ad{{\mathrm{ad}}}
\def\adth{\mathrm{ad}}

\def\Talpha#1{\vbox{\ialign{##\crcr
    $\alpha$\crcr\noalign{\kern2pt\nointerlineskip}
	   $\hfil\displaystyle{#1}\hfil$\crcr}}}

 \def\fin{\mathrm{End}}

\def\kb{{\mathbf{k}}}

\def\cala{{\mathcal A}}
\def\calb{{\mathcal B}}
\def\cald{{\mathcal D}}
\def\calk{{\mathcal K}}
\def\call{{\mathcal L}} 
\def\calm{{\mathcal M}} 

\def\calh{{\mathcal H}} 
\def\cals{{\mathcal S}} 
\def\calc{{\mathcal C}} 
\def\cale{{\mathcal E}}
\def\cals{{\mathcal S}}
\def\calu{{\mathcal U}}

\def\calf{{\mathcal F}} 
\def\calr{{\mathcal R}}
\def\calz{{\mathcal Z}}

\def\fracM{\mathfrak M} 
\def\fracf{\mathfrak f} 
\def\fracF{\mathfrak F}
 
\def\fracT{\mathfrak T}
\def\fracs{\mathfrak s} 
\def\fracS{\mathfrak S} 
\def\fracB{\mathfrak B}
\def\fraca{{\mathfrak A}} 
 
\def\fracg{{\mathfrak g}}
\def\fracl{{\mathfrak l}}
\def\fracu{{\mathfrak u}}
 
\def\bbbone{\mbox{\rm 1\hspace {-.6em} l}}

\def\der{\mathrm{\scriptsize Der}}
\def\derth{\mathrm{{\scriptsize Der}}}

\def\os{\underline{\Omega}}

\newtheorem{theo}{THEOREM}
\newtheorem{lemma}{LEMMA} 
\newtheorem{proposition}{PROPOSITION}

\begin{document}
    
    \begin{flushright}
	LPT-ORSAY 00/31
	\end{flushright}
	\vspace{1cm}

\title[Lectures on differentials, generalized differentials 
and...]{Lectures on differentials, generalized differentials \\ and 
on some examples related to theoretical physics}
\author{Michel Dubois-Violette}
\address{Laboratoire de Physique Th\'eorique\footnote{Unit\'e Mixte 
de Recherche du Centre National de la Recherche Scientifique - UMR 
8627}, B\^atiment 210, Universit\'e Paris XI, F-91405 Orsay Cedex}
\email{patricia@osiris.th.u-psud.fr}
\thanks{To be published in the Proceedings of Bariloche 2000, 
``Quantum Symmetries in Theoretical Physics and Mathematics", R. 
Coquereaux, R. Trinchero Eds, Contemporary Mathematics, American 
Mathematical Society.}

\begin{abstract}
These notes contain a survey of some aspects of the theory of 
differential modules and complexes as well as of their 
generalization, that is, the theory of $N$-differential modules and 
$N$-complexes. Several applications and examples coming from physics 
are discussed. The commun feature of these physical applications is 
that they deal with the theory of constrained or gauge systems. In 
particular different aspects of the BRS methods are explained and a 
detailed account of the $N$-complexes arising in the theory of higher 
spin gauge fields is given.
\end{abstract}
\maketitle

\section{Introduction}

Differential algebraic and (co)homological methods have rapidly
sprung up in theoretical physics in connection with the development 
of gauge theories. Their interventions occur at two levels, firstly 
at a classical level under a more systematic use of the calculus of 
differential forms, secondly under the emergence of the BRS methods 
in connection with the quantization of gauge theories. In fact the 
BRS technique provides an explicitely local and relativistic 
invariant way to develop perturbation theory for quantum gauge 
theories \cite{BRS}, \cite{BRS2}. It is worth noticing here that one 
cannot overestimate the role of the locality principle in 
perturbative renormalization \cite{EG}. Independently of these 
perturbative developments, methods for quantizing constrained systems 
on phase space have been developed using the path integral \cite{fv} 
which were obviously related. In both cases enter ``ghosts" \cite{FP} 
and the occurrence of a differential, i.e. an endomorphism of square 
zero. It turns out that the latter construction essentially reduces 
to a ``homological" description of classical constrained systems 
\cite{he:1} in which the ghosts and the differential have a natural 
interpretation in terms of standard mathematical concepts \cite{McM}, 
\cite{ros}, \cite{sta}, \cite{dv:2}.\\

Here we shall not give a systematic exposition of the above topics 
but, instead, we shall follow a sort of transversal way. These notes 
give a survey of appropriate concepts and results in homology which 
will be illustrated at each level with examples of application in 
theoretical physics. Furthermore recent developments in a 
generalization of homology will be reviewed as well as some  physical 
applications.\\

The plan is the following. In Section 2 we give the basic definitions 
and results on homology of differential modules. In Section 3 we 
introduce graduation, that is we discuss complexes and give several 
examples; in this section we explain the constructions connected with 
simplicial modules and we describe the tensor product of complexes. 
Section 4 is a physical illustration of the fact that there is no 
natural tensor product of differential modules whereas there is one 
for complexes; we show there that the introduction of ghosts at the 
one-particle level in the free field theory is worthwhile to render 
the theory natural over the physical space. In Section 5 we introduce 
$N$-differentials and discuss the generalization of homology 
associated with $N$-differential modules; we give there several 
examples of constructions some of which are related to physics (e.g. 
parafermions). Section 6 is devoted to the corresponding graded 
situation i.e. to $N$-complexes; we recall there the constructions of 
$N$-complexes associated to simplicial modules and the result which 
expresses in these cases the generalized homology in terms of the 
ordinary one (Theorem 2) \cite{D-V2}. In Section 7 which summarizes 
results of  \cite{D-VH}, \cite{D-VH2},  we introduce $N$-complexes of 
tensor fields on $\mathbb R^D$ generalizing the complex of 
differential forms and we state the corresponding generalization of 
the Poincar\'e lemma (Theorem 3); we also explain why these 
$N$-complexes naturally enter the theory of higher spin gauge fields. 
In Section 8 we discuss graded differential algebras and their 
``$N$-generalization" and give a universal $N$-construction 
generalizing the usual universal differential calculus over a unital 
associative algebra \cite{D-VK}, \cite{D-V2}. Section 9 describes the 
homological approach to ``subquotients" and applies it to constrained 
systems (BRS-method). The main result, Theorem 4, is slightly more 
general than the results of \cite{dv:2} (more general context), so we 
give a sketch of proof of Lemma~10 on which it relies.
Finally in Section 10 we generalize constructions of the previous 
section to $N$-differential modules in connection with a quantum 
gauge group problem arising for the zero modes in the 
Wess-Zumino-Novikov-Witten model; this section is a summary of 
\cite{D-VT2} (see also \cite{D-VT1}) .\\

These notes contain almost no proof, many results are classical or 
easy. There are two notable exceptions, namely Theorem 2 and Theorem 
3 the proof of which are absolutely non trivial although their 
meaning is transparent.\\

Let us say some words on our conventions. For sake of completeness we 
have given the formulation in terms of modules over a commutative 
ring $\kb$; the tensor product symbol $\otimes$ without other 
specification means the tensor product over $\kb$ (of $\kb$-modules), 
i.e. $\otimes=\otimes_\kb$. In the physical examples $\kb$ is either 
the field $\mathbb R$ of the field $\mathbb C$, so the reader may 
well understand $\kb$ like that and then the $\kb$-modules are vector 
spaces over $\mathbb R$ or $\mathbb C$. We use throughout the 
Einstein convention of summation of repeated up-down indices. A 
diagram of mappings between sets is said to be  a {\sl commutative 
diagram} if given two path of mappings between (two vertex) two sets 
of the diagram, the corresponding compositions of mappings coincide. 
A {\sl Young diagram} is a finite collection of boxes, or cells, 
arranged in left-justified rows, with a weakly decreasing number of 
cells in each row. Given a Young diagram of $n$ cells $Y$, one 
associates to it a projector $\mathbf{Y}$, the Young symmetrizer,  on 
the space of covariant tensors of degree $n$ on $\mathbb R^D$ by the 
following procedure. Let $T_{\mu_1\cdots\mu_n}$ be the components of 
$T$, then the components $\mathbf{Y}(T)_{\mu_1\cdots\mu_n}$ of 
$\mathbf{Y}(T)$ are obtained by filling successively the cells of the 
rows of $Y$ with $\mu_1,\cdots,\mu_n$, then by symmetrizing the 
$\mu$'s which belong to the same rows and then by antisymmetrizing 
the $\mu$'s which belong to the same columns. For Young diagrams 
etc., we use the notations of \cite{Ful}. We also mention that many 
subjects of these lectures are also treated in \cite{D-V4} so, 
although the aims of \cite{D-V4} are different, it is a complement 
for the present notes.

\section{Differential modules}

Throughout these notes, $\kb$ is a commutative ring with a unit and 
by a module without other specification, we always mean a 
$\kb$-module; the same convention is adopted for homomorphisms, 
endomorphisms, etc.. A module $E$ equipped with an endomorphism $d$ 
satisfying $d^2=0$ will be referred to as a {\sl differential module} 
and the endomorphism $d$ as its {\sl differential}. Given two 
differential modules $(E,d)$ and $(E',d')$, a {\sl homomorphism of 
differential modules} of $E$ into $E'$ is a homomorphism (of 
$\kb$-modules) $\varphi:E\rightarrow E'$ satisfying $\varphi\circ 
d=d'\circ \varphi$.\\

A sequence of homomorphisms of modules (resp. of differential modules)
\[
\cdots \longrightarrow E_i\stackrel{\varphi_i}{\longrightarrow} 
E_{i+1}\stackrel{\varphi_{i+1}}{\longrightarrow} 
E_{i+2}\longrightarrow \cdots
\]
is said to be {\sl exact} if $\im(\varphi_i)=\ker(\varphi_{i+1})$. In 
particular the sequence $0\rightarrow E 
\stackrel{\varphi}{\rightarrow} F$ is exact if and only if $\varphi$ 
is injective and the sequence $E\stackrel{\varphi}{\rightarrow} 
F\rightarrow 0$ is exact if and only if $\varphi$ is surjective.\\

Let $E$ be a differential module with differential $d$, then by 
definition one has $\im(d)\subset \ker(d)$ so the non exactness of 
the sequence $E\stackrel{d}{\rightarrow}E\stackrel{d}{\rightarrow}E$ 
is measured by the module $H(E)=\ker(d)/\im(d)$ which is referred to 
as the {\sl homology} of the differential module $E$. Let 
$\varphi:E\rightarrow F$ be a homomorphism of differential modules, 
then one has $\varphi(\im(d))\subset \im(d)$ and 
$\varphi(\ker(d))\subset \ker(d)$ (with an obvious abuse of 
notations) so $\varphi$ induces a homomorphism 
$\varphi_\ast:H(E)\rightarrow H(F)$ in homology. An important result 
for the computations of homology is given by the following 
proposition.

\begin{proposition}
Let $0\rightarrow E\stackrel{\varphi}{\rightarrow} F 
\stackrel{\psi}{\rightarrow} G\rightarrow 0$ be an exact sequence of 
differential modules; then there is a homomorphism 
$\partial:H(G)\rightarrow H(E)$ such that the triangle of 
homomorphisms

\[ \begin{diagram} 
\node{}\node{H(F)}\arrow{se,t}{\psi_\ast}\node{}\\
\node{H(E)}\arrow{ne,t}{\varphi_\ast}\node{}\node{H(G)}\arrow[2]{w,t}{\partial}
\end{diagram}
\]
is exact.
\end{proposition}

The exactness at $H(F)$ is easy and we only sketch the construction 
of $\partial$. Let $z\in G$ be such that $dz=0$ and let us denote by 
$[z]\in H(G)$ the class of $z$. Since $\psi$ is surjective there is a 
$y\in F$ such that $\psi(y)=z$; one has $\psi(dy)=d\psi(y)=dz=0$ so 
$dy\in\ker(\psi)$. By exactness at $F$, there is an $x\in E$ such 
that $\varphi(x)=dy$ and one has $\varphi(dx)=d\varphi(x)=d^2y=0$. 
Since $\varphi$ is injective it follows that $dx=0$ and we denote by 
$[x]\in H(E)$ the class of $x$. It turns out (and this is not 
difficult to verify) that $[x]\in H(E)$ does only depend on $[z]\in 
H(G)$ and that the mapping $[z]\mapsto \partial[z]=[x]$ is a 
homomorphism $\partial:H(G)\rightarrow H(E)$ which satisfies the 
statement of the proposition.\\

Quite generally, a five terms exact sequence of the form
\[
0\longrightarrow E\stackrel{\varphi}{\longrightarrow} 
F\stackrel{\psi}{\longrightarrow} G\longrightarrow 0
\]
is called a {\sl short exact sequence} and given a short exact 
sequence of differential modules as in Proposition 1, the 
homomorphism $\partial:H(G)\rightarrow H(E)$ is called the {\sl 
connecting homomorphism} of the short exact sequence of differential 
modules. The connecting homomorphism is natural (i.e. functorial) in 
the following sense: For any commutative diagram of differential 
modules
\[
\begin{diagram}
\node{0} \arrow{e} 
\node{E} \arrow{e,t}{\varphi}\arrow{s,l}{\lambda}
\node{F}\arrow{e,t}{\psi}\arrow{s,l}{\mu}
\node{G}\arrow{s,l}{\nu}\arrow{e}\node{0}\\
\node{0} \arrow{e} \node{E'} \arrow{e,t}{\varphi'} 
\node{F'} \arrow{e,t}{\psi'}
\node{G'}\arrow{e}\node{0}
\end{diagram}
\]
with exact rows, the diagram
\[
\begin{diagram}
\node{H(G)}\arrow{s,l}{\nu_\ast} \arrow{e,t}{\partial} 
\node{H(E)}\arrow{s,l}{\lambda_\ast}\\
\node{H(G')}\arrow{e,t}{\partial}\node{H(E')}
\end{diagram}
\]
is commutative.\\

It is worth noticing here that although direct sums of differential 
modules are well defined, there is no natural tensor product of 
differential modules. A natural tensor product will be only obtained 
in the graded case, that is for complexes (see below).\\

In the case where $\kb$ is a field, a differential module will be 
called a {\sl differential vector space} or simply a {\sl 
differential space}. In the examples connected with physics, $\kb$ 
will always be  either the field $\mathbb R$ or the field $\mathbb C$.

\section{Complexes}

By a {\sl complex}, without other specification, we always mean a 
differential module $E$ which is $\mathbb Z$-graded, 
$E=\displaystyle{\oplusinf_{n\in \mathbb Z}}E^n$, with a differential 
$d$ which is of degree 1 or $-1$. When $d$ is of degree $-1$, $E$ is 
referred to as a  {\sl chain complex} and when $d$ is of degree 1, 
$E$ is referred to as a {\sl cochain complex}. One passes from the 
chain complexes to the cochain ones by changing the signs of the 
degrees $(n\mapsto -n)$. In the following we shall only consider the 
cochain case. The homology of a cochain complex $E$ is usually 
referred to as the {\sl cohomology} of $E$. Since $d$ is homogeneous, 
the homology of a complex $E$ is $\mathbb Z$-graded~: 
$H(E)=\displaystyle{\oplusinf_{n\in \mathbb Z}}H^n(E)$ with 
$H^n(E)=\ker(d)\cap  E^n/\im(d)\cap E^n$. Many notions for complexes 
do only depend on the underlying $\mathbb Z_2$ graduation $(\mathbb 
Z_2=\mathbb Z/2\mathbb Z)$ so let us define a $\mathbb Z_2$-{\sl 
complex} to be a differential module $E$  which is $\mathbb 
Z_2$-graded, $E=E^0\oplus E^1$, with a differential $d$ which is of 
degree 1 (=$-1$). Again, the homology $H(E)$ of a $\mathbb 
Z_2$-complex is $\mathbb Z_2$-graded, that is $H(E)=H^0(E)\oplus 
H^1(E)$. A {\sl homomorphism of complexes} or 
{\sl of $\mathbb Z_2$-complexes} is a homomorphism of differential 
modules which is homogeneous of degree 0.\\

Let $0\longrightarrow E \stackrel{\varphi}{\longrightarrow} F 
\stackrel{\psi}{\longrightarrow} G\longrightarrow 0$ be a short exact 
sequence of cochain complexes; it follows from the definition of the 
connecting homomorphism $\partial$ that the exact triangle of 
Proposition 1 gives rise to the long exact sequence of homomorphisms
\[
\cdots \stackrel{\partial}{\longrightarrow} H^n(E) 
\stackrel{\varphi_\ast}{\longrightarrow} H^n(F) 
\stackrel{\psi_\ast}{\longrightarrow} H^n(G) 
\stackrel{\partial}{\longrightarrow}H^{n+1}(E)\stackrel{\varphi_\ast}{\longrightarrow}\cdots
\]
in cohomology. Similarily if $0\longrightarrow E 
\stackrel{\varphi}{\longrightarrow} F\stackrel{\psi}{\longrightarrow} 
G\longrightarrow 0$ is a short exact sequence of $\mathbb 
Z_2$-complexes, the exact triangle of Proposition 1 gives rise to the 
exact hexagon of homomorphisms
\[ 
\begin{diagram} 
\node{}
\node{H^0(F)} \arrow{e,t}{\psi_\ast}\node{H^0(G)}
\arrow{se,t}{\partial} \node{} \\ 
\node{H^0(E)}
\arrow{ne,t}{\varphi_\ast}
\node{} \node{} \node{H^1(E)} \arrow{sw,b}{\varphi_\ast} \\ \node{}
\node[1]{H^1(G)} \arrow{nw,b}{\partial}
\node{H^1(F)}
\arrow{w,b}{\psi_\ast} \node{} \end{diagram} 
\]
for the corresponding homologies.\\

Let $E$ and $F$ be two cochain complexes, (resp. $\mathbb 
Z_2$-complexes), {\sl their tensor product} $E\otimes F$ is the 
graded module $E\otimes F=\displaystyle{\oplusinf_n}(E\otimes F)^n$ 
with $(E\otimes F)^n=\displaystyle{\oplusinf_{r+s=n}}E^r\otimes F^s$ 
equipped with the differential $d$ defined by
\[
d(e\otimes f) =de \otimes f + (-1)^n e\otimes df,
\]
for any $e\in E^n$ and $f\in F$. One verifies that so defined on 
$E\otimes F$, $d$ is homogeneous of degree 1 and satisfies $d^2=0$ so 
that $E\otimes F$ is again a cochain complex, (resp. a $\mathbb 
Z_2$-complex). The virtue of this definition is the K\"unneth formula 
which we describe only for complexes of vector spaces in the 
following proposition, \cite{ghv}, \cite{weib}.

\begin{proposition}
Assume that the ring $\kb$ is a field then one has $H(E\otimes 
F)=H(E)\otimes H(F)$.
\end{proposition}

The above tensor product being the tensor product of graded vector 
spaces (over $\kb$) i.e. $H^n(E\otimes 
F)=\displaystyle{\oplusinf_{r+s=n}}H^r(E)\otimes H^s(F)$. This 
formula applies as well to the (co)chain complexes case and to the 
$\mathbb Z_2$-complexes case (whenever $\kb$ is a field).\\

In the next section we shall describe a physical application of 
Proposition 2 combined with the remark that there is no such tensor 
product for differential spaces. We now achieve this section by the 
description of some classical constructions which will be used 
later.\\

Let $\fracg$ be a Lie algebra, let $R$ be a representation 
space  of  $\fracg$ and denote by  $X\mapsto \pi(X)\in \fin(R)$
the action of $\fracg$ on $R$. An $R$-{\sl  valued (Lie
algebra) $n$-cochain of $\fracg$} is a linear mapping 
$X_1\wedge\dots\wedge X_n\mapsto \omega(X_1,\dots,X_n)$ of
$\bigwedge^n\fracg$ into $R$.  The vector space of these
$n$-cochains will be denoted by  $C^n_\wedge(\fracg,R)$. One
defines a homogeneous endomorphism $d$ of degree 1 of the
$\mathbb  N$-graded vector space
$C_\wedge(\fracg,R)=\oplusinf_n C^n_\wedge(\fracg,R)$ of  all
$R$-valued cochains of $\fracg$ by setting \[ \begin{array}{ll}
d(\omega)(X_0,\dots,X_n) & =\sum^n_{k=0}(-1)^
k\pi(X_k)\omega(X_0,\stackrel{k\atop \vee}{\dots},X_n)\\ & 
+\sum_{0\leq r<s\leq
n}(-1)^{r+s}\omega([X_r,X_s],X_0\stackrel{r\atop\vee}{\dots} 
\stackrel{s\atop \vee}{\dots}  X_n) \end{array} \] for
$\omega\in  C^n_\wedge(\fracg,R)$ and $X_i\in\fracg$. It
follows from the Jacobi identity and from
$\pi(X)\pi(Y)-\pi(Y)\pi(X)=\pi([X,Y])$ that $d^2=0$. Thus
equipped  with $d$, $C_\wedge(\fracg,R)$ is a cochain complex
and its cohomology, denoted  by $H(\fracg,R)$, is called the
$R$-{\sl valued cohomology of $\fracg$}. The complexes 
$C_\wedge(\fracg,R)$ are also called {\sl Chevalley-Eilenberg 
complexes} and the differential $d$ is the {\sl Chevalley-Eilenberg 
differential}.\\

There is a standard way to produce positive complexes (i.e. complexes 
$E=\oplus E^n$ with $E^n=0$ for $n<0$) 
starting from (co)simplicial modules, (see e.g. \cite{jll}, 
\cite{weib}). A
{\sl pre-cosimplicial module} (or {\sl semi-cosimplicial} in the 
terminology of
\cite{weib}) is a sequence of modules $(E^n)_{n\in \mathbb N}$ 
together with
{\sl coface homomorphisms} ${\mathfrak f}_i: E^n\rightarrow E^{n+1},\ 
i\in
\{0,1,\dots,n+1\}$, satisfying\\
$({\fracF})$\hspace{2cm} ${\fracf}_j 
{\fracf}_i={\fracf}_i{\fracf}_{j-1}$\hspace{.5cm} if
$i<j$.\\
Given a pre-cosimplicial module $(E^n)$, one associates to it a 
positive
complex $(E,d)$ by setting $E=\oplus_{n\in \mathbb N} E^n$ and $
d=\sum^{n+1}_{i=0} (-1)^i \fracf_i: E^n\rightarrow E^{n+1}$.
One verifies that $d^2=0$ is implied by the coface relations 
$(\fracF)$. The differential $d$ will be referred to as the {\sl 
simplicial differential} of $(E^n)$.
The cohomology $H(E)=\oplus H^n(E)$ with 
$H^n(E)=\ker(d:E^n\rightarrow E^{n+1})/dE^{n-1}$
 of $(E,d)$  will be referred to as {\sl the
cohomology of the pre-cosimplicial module} $(E^n)$.
A {\sl cosimplicial
module} is a pre-cosimplicial module $(E^n)$ with coface 
homomorphisms $\fracf_i$
as before together with {\sl codegeneracy homomorphisms}
$\fracs_i:E^{n+1}\rightarrow E^n,\ i\in \{0,\dots,n\}$, satisfying\\
$(\fracS)\hspace{1.5cm} \fracs_j\fracs_i=\fracs_i 
\fracs_{j+1}\hspace{1cm}$ if $i\leq j$\\
and\\
$(\fracS\fracF)\hspace{1.5cm} \fracs_j\fracf_i=\left\{
\begin{array}{lll}
\fracf_i\fracs_{j-1} & \mbox{if} & i<j\\
I & \mbox{if} & i=j\  \mbox{or}\  i=j+1\\
\fracf_{i-1}\fracs_j & \mbox{if} & i>j+1
\end{array}
\right.
$\\
Given a cosimplicial module $(E^n)$ the elements $\omega$ of $E^n$ 
such that $s_i(\omega)=0$ for $i\in \{0,\cdots,n\}$ are called 
normalized cochains of degree $n$  and the graded module 
$N(E)=\displaystyle{\oplusinf_n} N^n(E)$ of all normalized cochains 
is a subcomplex of $E$ which has the same cohomology as the one of 
$E$, i.e. $H(E)$. The correspondence $(E^n)\mapsto N(E)$ defines an 
equivalence between the category of cosimplicial modules and the 
category of positive cochain complexes \cite{weib} which is referred 
to as {\sl the Dold-Kan correspondence} (for the category of 
$\kb$-modules).\\

Let $\cala$ be an associative unital $\kb$-algebra and let $\calm$ be 
an
$(\cala,\cala)$-bimodule. A $\calm$-{\sl valued  Hochschild cochain}
of degree $n$ or {\sl Hochschild $n$-cochain} of $\cala$ is a linear 
mapping 
$x_1\otimes \cdots \otimes x_n\mapsto \omega(x_1,\cdots, x_n)$
of $\otimes^n\cala$ into $\calm$. The $\kb$-module of
all $\calm$-valued Hochschild $n$-cochains is denoted by 
$C^n(\cala,\calm)$. The
sequence $(C^n(\cala,\calm))$ is a cosimplicial module with cofaces 
$\fracf_i$ and
codegeneracies $\fracs_i$ defined by \cite{jll}, \cite{weib}\\
$\fracf_0(\omega)(x_0,\dots,x_n)=x_0\omega(x_1,\dots,x_n)$\\
$\fracf_i(\omega)(x_0,\dots,x_n)=\omega(x_0,\dots,x_{i-1}x_i,\dots,x_n)$\hspace{2cm} 
for $i\in\{1,\dots,n\}$\\
$\fracf_{n+1}(\omega)(x_0,\dots,x_n)=\omega(x_0,\dots,x_{n-1})x_n$\\
and\\
$\fracs_i(\omega)(x_1,\dots,x_{n-1})=\omega(x_1,\dots,x_i,\bbbone,x_{i+1},\dots,x_{n-1})$ 
\hspace{1cm} for $i\in \{0,\dots,n-1\}$\\
for $\omega\in C^n(\cala,\calm)$ and $x_i\in \cala$. The cohomology 
$H(\cala,\calm)$ of this cosimplicial module is the {\sl 
$\calm$-valued Hochschild cohomology of $\cala$}. In his case the 
simplicial differential is called the {\sl Hochschild differential}.\\

There is a relation between the cohomology of a Lie algebra $\fracg$ 
and the Hochschild cohomology of its universal enveloping algebra 
$U(\fracg)$  which we now describe again in the case where $\kb$ is a 
field. Given a bimodule $\calm$ over $U(\fracg)$ (that is 
a\linebreak[4]  $(U(\fracg),U(\fracg))$-bimodule), let us define the 
representation $X\mapsto \ad (X)$ of $\fracg$ in the vector space 
$\calm$ by $\ad(X)m=Xm-mX$ for $X\in\fracg$ and $m\in\calm$. Let 
$H(\fracg,\calm^{\adth})$ denote the Lie algebra cohomology of 
$\fracg$ with values in $\calm$ for the $\ad$ representation; its 
relation with the $\calm$-valued Hochschild cohomology of 
$U(\fracg)$, $H(U(\fracg),\calm)$ is given by the following theorem 
\cite{cart:02}, \cite{jll}.

\begin{theo}
Assume that $\kb$ is a field, let $\fracg$ be a Lie algebra over 
$\kb$ and let $\calm$ be a bimodule over $U(\fracg)$. Then there is a 
canonical isomorphism $H(\fracg,\calm^{\adth})\simeq 
H(U(\fracg),\calm)$.
\end{theo}

If $R$ is a representation space of $\fracg$ with action $X\mapsto 
\pi(X)$, then by the very definition of $U(\fracg)$,  $\pi$ extends 
as a representation of $U(\fracg)$ so $R$ is canonically  a left 
$U(\fracg)$-module. One converts $R$ into a 
$(U(\fracg),U(\fracg))$-bimodule $\calr$ by acting on the right with 
the trivial action given by the counit of $U(\fracg)$ (recall that 
$U(\fracg)$ is a Hopf algebra); one then has $R=\calr^{\adth}$.

\section{A physical example: Naturality of ghosts}

The Wigner one-particle space for mass zero and spin one is the
direct hilbertian integral $\int_{C_+} d\mu_0(p)\calh(p)$ of
2-dimensional Hilbert spaces $\calh(p)$ over the future light
cone \[ C_+=\{p\vert g^{\mu\nu}p_\mu p_\nu=p_0^2-\vec p^2=0,\ \
p_0>0\} \] with respect to the invariant measure
$d\mu_0(p)=\frac{1}{(2\pi)^3}\frac{d^3\vec p}{2p^0}$, where
$\calh(p)$ is the quotient of the subspace $\calz(p)=\{A_\mu\in
\mathbb C^4\vert p^\mu A_\mu=0\}$ of $\calc(p)=\mathbb C^4$ by
the subspace  $\calb(p)=\{p_\mu\varphi\vert \varphi\in \mathbb
C\}$ spanned by $p$, the scalar product of $\calh(p)$ being
induced by the indefinite scalar product of $\calc(p)$ defined
by $\langle A\vert A'\rangle=-g^{\mu\nu}\bar A_\mu A'_\nu$. The
scalar product of $\calc(p)$ is positive semi-definite on
$\calz(p)$ and $\calb(p)$ is its isotropic subspace whereas
$\calz(p)$ is the orthogonal of $\calb(p)$ in $\calc(p)$. Notice
that the indefinite metric space $\calc(p)$ does not depend on
$p$; we keep the reference to $p$ in order to remember that it
carries a representation of {\sl the little group at} $p$.  The 
little group at $p$ here means the subgroup $\call_p$ of the Lorentz 
group which consists of the Lorentz tranformations $\Lambda$ 
preserving the (quadri) vector $p$, that is
\[
\call_p=\{\Lambda\in GL(4,\mathbb R)\ \ \vert\ \  \Lambda^\mu_\lambda 
\Lambda^\nu_\rho g^{\lambda\rho}=g^{\mu\nu}\ \ \mbox{and}\ \ 
\Lambda^\mu_\nu p^\nu=p^\mu\}.
\]
The occurrence of such a triplet $(\calc(p), \calz(p),\calb(p))$
where $\calc(p)$ has an indefinite scalar product with
$\calb(p)$ isotropic having $\calz(p)$ as orthogonal, etc. is
familiar in connection with indecomposable representations of
groups (here the little group) \cite{G.R}, \cite{H.A} and the
indefinite metric is furthermore required to get a local
covariant description of the electromagnetic gauge potential
\cite{Stroc}, \cite{SW}, see also \cite{J-PL} in this context.

Let $Q(p)=Q$ be the linear endomorphism of $\calc(p)$ defined by
$Q(A)_\mu=p_\mu p^\nu A_\nu$. Then $Q$ is hermitian, i.e.
$\langle A\vert QA'\rangle = \langle Q A\vert A'\rangle$, and
one has $Q^2=0$ in view of $p_\mu p^\mu=0$. Furthermore the
image of $Q$ is $\calb(p)$ and its kernel is $\calz(p)$. In
other words $(\calc(p),Q(p)$) is a differential space
 and $\calh(p)$ is its homology, i.e. one has
$\calh(p)=\ker(Q)/\im(Q)$. Thus, apart from questions of domain
and function spaces, everything is perfect at the
``one-particle" level: Namely one has an indefinite metric space
$\calc$ which consists of functions $p\mapsto A_\mu(p)\in
\calc(p)$ on the light cone $C_+$ and which is equipped with a
differential $Q$ (i.e. an endomorphism satisfying $Q^2=0$) such
that the physical one-particle space, (i.e. the Wigner space),
is the homology $\ker(Q)/\im(Q)$ of $\calc$.\\

As is well known, the role of $\calc$ is to provide, via the Fock
space constructions, an indefinite metric space on which the
local covariant gauge potential (free) field operator acts; the
corresponding space of physical states being of course the Fock
space constructed over the one-particle Wigner space. However it
turns out that the above one-particle (homological) picture does
not generalize naively at the $n$-particle level for $n\geq 2$.
To show what is involved here, let us analyze the situation at
the two-particle level. In order to avoid complications
connected with the problem of the choice of the function space
and with the problem of symmetrization, let us work at fixed
momenta $p_1$ and $p_2$ on the light cone $C_+$ with $p_1\not=
p_2$. The indefinite metric space is then the 16-dimensional
space $\calc(p_1)\otimes \calc(p_2)$ whereas the space of
physical states is the 4-dimensional Hilbert space
$\calh(p_1)\otimes\calh (p_2)$. The point now is that there is
no canonical way to construct $\calh(p_1)\otimes \calh(p_2)$
from $\calc(p_1)\otimes \calc(p_2)$. More precisely there is no
canonical way to build a differential on 
$\calc(p_1)\otimes\calc(p_2)$ out
of the differentials $Q(p_1)$ and $Q(p_2)$ of $\calc(p_1)$ and
$\calc(p_2)$ in such a way that its homology is
$\calh(p_1)\otimes \calh(p_2)$. In fact the most natural
candidate would be $Q_{12}=Q(p_1)\otimes \id_2+\id_1\otimes
Q(p_2)$ but this is not of square zero, only its third power
vanishes, (for the ``n-particle" case it would be the $(n+1)$-th
power). Notice that with $Q_{12}$ satisfying $(Q_{12})^3=0$ one
can associate the generalized homologies (see below)
$H_{(1)}(Q_{12})=\ker(Q_{12})/\im((Q_{12})^2)$ and
$H_{(2)}(Q_{12})=\ker((Q_{12})^2)/\im(Q_{12})$ however it is
easy to show that one canonically has
$H_{(1)}(Q_{12})=\calz(p_1)\otimes \calz(p_2)$ and that
$H_{(2)}(Q_{12})$ is isomorphic to $H_{(1)}(Q_{12})$. Thus
$H_{(1)}(Q_{12})$ is a subspace of $\calc(p_1)\otimes \calc(p_2)$ on
which the metric is positive semi-definite but it is still not
the physical space $\calh(p_1)\otimes \calh(p_2)$.\\

Notice that we do not claim that there is no differential on 
$\calc(p_1)\otimes \calc(p_2)$ such that the corresponding homology 
is $\calh(p_1)\otimes \calh(p_2)$ but that we claim that there is no 
canonical one, that is no reasonable expression for such a 
differential in terms of the differentials $Q(p_1)$ and $Q(p_2)$. We 
refer to Appendix A for the precise statement.\\

As pointed out above, the origin of the difficulty is the
non-existence of a good tensor product between differential
spaces, i.e. between vector spaces equipped with endomorphisms
of square zero. If instead of differential spaces one has
complexes (of vector spaces), then the
situation is much better; namely one has a canonical tensor
product of complexes which is such that the homology of the
tensor product is the tensor product of the homologies, (see
last section). Furthermore one can show 
that the symmetrization-antisymmetrization involved in the Fock
space construction does not spoil this picture.\\

Fortunately there is a canonical way (related to Theorem 4) to 
construct a complex
$C(p)=C^{-1}(p)\oplus C^0(p)\oplus C^1(p)$ with a differential
of degree 1 such that $C^0(p)=\calc(p)$ and such that its
(co)homology is again $\calh(p)$. We now describe this construction. 
Let
$\varepsilon^\mu$ be the (real) canonical base of
$\calc(p)=C^0(p)=\mathbb C^4$ and let $\omega^{(+)}$ and
$\omega^{(-)}$ be the basis of the one dimensional spaces
$C^1(p)$ and $C^{-1}(p)$ $(\cong \mathbb C)$. Define the
homogeneous linear endomorphism $\delta(p)=\delta$ of degree 1
of $C(p)$ by $\delta\omega^{(+)}=0$,
$\delta\varepsilon^\mu=\alpha p^\mu\omega^{(+)}$ and
$\delta\omega^{(-)}=p_\mu\varepsilon^\mu$, ($\alpha$ being a
non-vanishing constant). It is clear that $\delta^2=0$ and it is
straightforward to verify that the (co)homology
$H(C(p))=\ker(\delta)/\im(\delta)$ of $C(p)$ is given by
$H(C(p))=H^0(C(p))=\calh(p)$. Notice that if
$c\omega^{(-)}+A_\mu\varepsilon^\mu+\tilde c\omega^{(+)}$ is an
arbitrary element of $C(p)$, $\delta$ reads in components
$\delta A_\mu=p_\mu c$, $\delta c=0$ and $\delta\tilde c=\alpha
p^\lambda A_\lambda$. One defines an indefinite hermitian
scalar product on $C(p)$ extending the one of $C^0(p)=\calc(p)$
for which $\delta$ is hermitian by setting $\langle
\varepsilon^\mu\vert \varepsilon^\nu\rangle=-g^{\mu\nu},$
$\langle\omega^{(+)}\vert \varepsilon^\mu\rangle=0,$
$\langle\omega^{(+)}\vert \omega^{(+)}\rangle=0,$
$\langle\omega^{(-)}\vert \varepsilon^\mu\rangle=0$,
$\langle\omega^{(-)}\vert \omega^{(-)}\rangle=0$ and
$\langle\omega^{(-)}\vert \omega^{(+)}\rangle=-\alpha^{-1}$.
 One can now
construct the generalized Fock space $\fracF(C)$ over the graded space
$C$ of ``functions" $p\mapsto (\tilde c(p), A_\mu(p), c(p))\in
C(p)$ on the future light cone. The space $\fracF(C)$ is the
graded-commutative algebra (freely) generated by the graded
vector space $C$ and one extends $\delta$ as an antiderivation
of $\fracF(C)$, again denoted by $\delta$, which still satisfies
$\delta^2=0$. The scalar product of $C$ extends canonically into an
indefinite scalar product of $\fracF(C)$ for which $\delta$ is
hermitian and the cohomology $H^0(\delta)$ is (a dense subspace
of) the physical space (i.e. the Fock space over the Wigner
one-particle space). One then constructs as usual the local
gauge potential field operator corresponding to the above
one-particle $A_\mu$ as well as the fermionic ghost and
antighost field operators corresponding to the above
one-particle $c$ and $\tilde c$. In order that the ghost and the
antighost fields be relatively local, it is necessary to take
$\alpha$ purely imaginary, i.e. $\alpha=i\lambda$ with
$\lambda\in \mathbb R_\ast$, otherwise one would obtain a factor
$D^{(1)}$ in their anticommutators. With this choice
($\alpha=i\lambda, \lambda\in \mathbb R_\ast$) the gauge
potential, the ghost and the antighost field operators are local
and relatively local, (see e.g. in \cite{NO}). Moreover these
fields are hermitian by their very definition.\\

Let us say a few words on the
case of spin two (and zero mass). In this case, the Wigner
one-particle space is again the direct hilbertian integral
$\int_{C_+}d\mu_0(p)\calh(p)$ of two-dimensional Hilbert
spaces $\calh(p)$ over the future light cone with respect to
$d\mu_0$ with $\calh(p)=\calz(p)/\calb(p)$ and $\calz(p)
\subset \calc(p)$ as above but now, $\calc(p)$ is the \linebreak[4]
10-dimensional space of symmetric tensors $h_{\mu\nu}=h_{\nu\mu}$,
\[
\begin{array}{lll}
\calz(p) & = & \{h_{\mu\nu}\in \calc(p)
\vert p^\mu(h_{\mu\nu}-\frac{1}{2}g_{\mu\nu}
g^{\alpha\beta}h_{\alpha\beta})=0\},\\
\\
\calb(p) & = &
\{p_\mu\varphi_\nu+p_\nu\varphi_\mu\vert \varphi_\lambda\in
\mathbb C^4\}
\end{array}\
\]
and the scalar product of $\calh(p)$ is induced
by the indefinite scalar product of $\calc(p)$ defined by
$\langle h\vert h'\rangle=g^{\mu\nu}g^{\lambda \rho}\bar
h_{\mu\lambda} h'_{\nu\rho}-\frac{1}{2}g^{\alpha\beta}\bar
h_{\alpha\beta} g^{\gamma\delta}h'_{\gamma \delta}$. Again
$\calb(p)$ is a completely isotropic (4-dimensional) subspace
of $\calc(p)$ whereas the 6-dimensional space $\calz(p)$ is
its orthogonal in $\calc(p)$, $(\calz(p)=\calb(p)^\perp)$. It
is worth noticing here that, apart from a multiplicative
constant, the scalar product $\langle h\vert h'\rangle$ is the
unique non-trivial covariant scalar product on $\calc(p)$ for
which $\calb(p)$ is isotropic; equivalently, the condition
$p^\mu(h_{\mu\nu}-\frac{1}{2}
g_{\mu\nu}g^{\alpha\beta}h_{\alpha \beta})=0$ is the unique
covariant linear (gauge) condition preserved by the
translations of $\calb(p)$. In view of the connection between
the classical linearized gravity theory and the massless spin
two particle, it is natural to interpret $h_{\mu\nu}\in
\calc(p)$ as the positive frequency part of the Fourier
transform at $p$ of a perturbation
$\underline{g}_{\mu\nu}(x)=g_{\mu\nu}+\varepsilon
h_{\mu\nu}(x)$ of the Minkowskian metric $g_{\mu\nu}$.
Translations by $\calb(p)$ then read $h_{\mu\nu}(x)\mapsto
h_{\mu\nu}(x)+\partial_\mu\varphi_\nu(x)+\partial_\nu\varphi_\mu(x)$
which corresponds to the first order in  $\varepsilon$ (i.e.
the linearization) of the action of infinitesimal
diffeomorphisms (i.e. vector fields) whereas the condition to
be in $\calz(p)$ translates into
$\partial^\mu(h_{\mu\nu}(x)-\frac{1}{2} g_{\mu\nu}
g^{\alpha\beta} h_{\alpha\beta}(x))=0$ which is the first
order in $\varepsilon$ of the de Donder harmonic coordinates
condition
$\frac{1}{\sqrt{\vert\underline{g}\vert}}\partial_\mu
\left(\sqrt{\vert\underline
g\vert}\underline{g}^{\mu\nu}\right)=\Delta_{\underline{g}}(x^\nu)=0$.
It may well be that this observation (i.e. connection between
Poincar\'e covariant Wigner analysis and de Donder harmonic
coordinates condition) is a little more than a curiosity. In
any case, we can now proceed as for the spin one case. One
defines the graded vector space $C(p)=C^{-1}(p)\oplus
C^0(p)\oplus C^1(p)$ by $C^0(p)=\calc(p)$ and $C^{-1}(p)\simeq
\mathbb C^4\simeq C^1(p)$ and we let $\omega^{(-)\mu}$ and
$\omega^{(+)\mu}$ be the basis of $C^{-1}(p)$ and $C^{+1}(p)$
corresponding to the canonical base $\varepsilon^\mu$ of
$\mathbb C^4$ and
$\varepsilon^{\mu\nu}=\frac{1}{2}(\varepsilon^\mu\otimes
\varepsilon^\nu+\varepsilon^\nu\otimes \varepsilon^\mu)$ be
the corresponding basis of $C^0(p)=\calc(p)$. One defines then
a differential $\delta$ of degree 1 of $C(p)$ by setting
$\delta\omega^{(+)\mu}=0$,
$\delta\varepsilon^{\mu\nu}=\alpha(p^\mu\omega^{(+)\nu}+
p^\nu\omega^{(+)\mu})$ and
$\delta\omega^{(-)\mu}=p_\nu(\varepsilon^{\mu\nu}-\frac{1}{2}
g^{\mu\nu}g_{\alpha\beta}\varepsilon^{\alpha\beta})$,
$\alpha\in \mathbb C_\ast$. Again one verifies that the cohomology 
$H(C(p))=\ker(\delta)/\im (\delta)$ of $C(p)$ is given by 
$H(C(p))=H^0(C(p))=\calh(p)$. If we let
$c_\rho\omega^{(-)\rho}+h_{\mu\nu}\varepsilon^{\mu\nu}+\tilde
c_\lambda\omega^{(+)\lambda}$ be an arbitrary element of
$C(p),\delta$ reads in components $\delta
h_{\mu\nu}=p_\mu c_\nu+p_\nu c_\mu$, $\delta c_\mu=0$ and
$\delta \tilde c_\mu=\alpha
p^\nu(h_{\mu\nu}-\frac{1}{2}
g_{\mu\nu}g^{\alpha\beta}h_{\alpha\beta})$.
Finally, one defines an indefinite hermitian scalar product
on $C(p)$ extending the one of $C^0(p)=\calc(p)$ for which
$\delta$ is hermitian by setting
$\langle\varepsilon^{\lambda\rho}\vert\varepsilon^{\mu\nu}
\rangle=\frac{1}{2}(g^{\lambda\mu}g^{\rho
\nu}+g^{\lambda\nu}g^{\rho\mu}-g^{\lambda\rho}g^{\mu\nu})$,
$\langle\omega^{(+)\lambda}\vert
\varepsilon^{\mu\nu}\rangle=0$, $\langle\omega^{(+)\mu}\vert
\omega^{(+)\nu}\rangle=0$,
$\langle\omega^{(-)\lambda}\vert\varepsilon^{\mu\nu}\rangle=0$,
$\langle\omega^{(-)\mu}\vert \omega^{(-)\nu}\rangle=0$ and
$\langle\omega^{(-)\mu}\vert
\omega^{(+)\nu}\rangle=\frac{1}{2\alpha}g^{\mu\nu}$. Thus,
apart from numbers of components, everything works as in the
case of spin one, in particular one must again take
$\alpha=i\lambda$ with $\lambda\in \mathbb R_\ast$ in order to
have locality and relative locality between the hermitian free
fields corresponding to $h_{\mu\nu}, c_\lambda$ and $\tilde
c_\rho$.\\

The main message of this section is ``the natural
necessity" of ghosts (i.e. of graduations) in order to have a
canonical local construction over the physical space and the
fact that, in the previous examples (and others), there is a
canonical way to introduce their counterpart at the
one-particle level. This rewriting of the free field theory for
zero mass and spin $\geq 1$ is certainly needed in order to
start to introduce consistently interactions between abelian
gauge fields. In particular this reformulation can be
considered as the zero-step for the perturbative construction
of quantum operatorial Yang-Mills theory.

\section{$N$-differential modules}

In the following, $N$ is a positive integer with $N\geq 2$. A module 
$E$ equipped with an endomorphism $d$ satisfying $d^N=0$ will be 
referred to as an $N$-{\sl differential module} and the endomorphism 
$d$ as its $N$-{\sl differential}. With this terminology, a 
2-differential module is just a differential module.For
each integer $m$ with $1\leq m\leq N-1$, one defines the sub-modules
$Z_{(m)}(E)$ and $B_{(m)}(E)$ by setting $Z_{(m)}(E)=\ker (d^m)$ and
$B_{(m)}(E)=\im(d^{N-m})$. It follows from the equation $d^N=0$ that 
$B_{(m)}(E)$ is a
submodule of $Z_{(m)}(E)$ and the quotient modules
$H_{(m)}(E)=Z_{(m)}(E)/B_{(m)}(E)$,\linebreak[4] $m\in 
\{1,\dots,N-1\}$, will be referred to
as {\sl the (generalized) homology} of the\linebreak[4]  
$N$-differential module $E$.\\

Let $m$ be an integer with $1\leq m\leq N-2$ and let $E$ be an 
$N$-differential
module. One has the inclusions $Z_{(m)}(E)\subset Z_{(m+1)}(E)$ and 
$B_{(m)}(E)\subset
B_{(m+1)}(E)$ which induces a homomorphism $[i]:H_{(m)}(E)\rightarrow
H_{(m+1)}(E)$. One  has also the inclusions $dZ_{(m+1)}(E)\subset 
Z_{(m)}(E)$ and
$dB_{(m+1)}(E)\subset B_{(m)}(E)$ which induces a homomorphism
$[d]:H_{(m+1)}(E)\rightarrow H_{(m)}(E)$. The following basic result 
show that
the $H_{(m)}(E)$ are not independent \cite{D-VK}, \cite{D-V2}.

\begin{lemma}

Let $\ell$ and $m$ be integers with $\ell\geq 1,\ m\geq 1$ and 
$\ell+m\leq
N-1$. Then the following hexagon $(\calh^{\ell,m})$ of homomorphisms
\[ \begin{diagram} \node{}
\node{H_{(\ell+m)}(E)} \arrow{e,t}{[d]^m} \node{H_{(\ell)}(E)}
\arrow{se,t}{[i]^{N-(\ell+m)}} \node{} \\ \node{H_{(m)}(E)}
\arrow{ne,t}{[i]^\ell}
\node{} \node{} \node{H_{(N-m)}(E)} \arrow{sw,b}{[d]^\ell} \\ \node{}
\node[1]{H_{(N-\ell)}(E)} \arrow{nw,b}{[d]^{N-(\ell+m)}}
\node{H_{(N-(\ell+m))}(E)}
\arrow{w,b}{[i]^m} \node{} \end{diagram} \]
is exact.
\end{lemma}

One has obvious notions of homomorphism of $N$-differential modules,
of\linebreak[4] $N$-differential submodule of an $N$-differential 
module, etc..
Let $\varphi:E\rightarrow E'$ be a homomorphism of $N$-differential 
modules.
Then one has $\varphi(Z_{(m)}(E))\subset Z_{(m)}(E')$ and
$\varphi(B_{(m)}(E))\subset B_{(m)}(E')$ which implies that $\varphi$ 
induces
a homomorphism\linebreak[4]   $\varphi_\ast:H_{(m)}(E)\rightarrow
H_{(m)}(E')$, $\forall m\in \{1,\dots, N-1\}$. Moreover $\varphi_\ast$
satisfies $\varphi_\ast \circ [i]=[i]\circ \varphi_\ast$ and 
$\varphi_\ast\circ
[d]=[d]\circ \varphi_\ast$. Proposition 1 has the following 
generalization for $N$-differential modules.

\begin{proposition}
Let $0\rightarrow E
\stackrel{\varphi}{\rightarrow}F\stackrel{\psi}{\rightarrow} G 
\rightarrow 0$
be a short exact sequence of\linebreak[4] $N$-differential modules. 
Then there
are homomorphisms $\partial:H_{(m)}(G)\rightarrow H_{(N-m)}(E)$ for 
$m\in
\{1,\dots, N-1\}$ such that the following hexagons $(\calh_n)$ of 
homomorphisms

\[ \begin{diagram} \node{}
\node{H_{(n)}(F)} \arrow{e,t}{\psi_\ast} \node{H_{(n)}(G)}
\arrow{se,t}{\partial} \node{} \\ \node{H_{(n)}(E)} 
\arrow{ne,t}{\varphi_\ast}
\node{} \node{} \node{H_{(N-n)}(E)} \arrow{sw,b}{\varphi_\ast} \\ 
\node{}
\node[1]{H_{(N-n)}(G)} \arrow{nw,b}{\partial} \node{H_{(N-n)}(F)}
\arrow{w,b}{\psi_\ast} \node{} \end{diagram} \]

are exact, for $n\in \{1,\dots,N-1\}$.
\end{proposition}

For a proof, we refer to \cite{KW}, \cite{D-V2}, \cite{D-V3}. In 
fact, there is a way to interpret $(\calh_n)$ as the exact hexagon 
corresponding to a short exact sequence of $\mathbb Z_2$-complexes 
$0\rightarrow C_{(n)}(E)\rightarrow C_{(n)}(F) \rightarrow 
C_{(n)}(G)\rightarrow 0$  associated with the $N$-complexes, 
\cite{D-V3}.\\

Let us now give  some criteria ensuring the vanishing of the 
$H_{(n)}(E)$. The
first criterion is extracted from \cite{Kap}.

\begin{lemma}
Let $E$ be an $N$-differential module such that $H_{(k)}(E)=0$ for 
some integer
$k$ with $1\leq k\leq N-1$. Then one has $H_{(n)}(E)=0$ for any 
integer $n$
with $1\leq n\leq N-1$.
\end{lemma}

A short proof of this lemma using Lemma 1 is given in \cite{D-V2}.
The next criterion which is in \cite{Kap} is connected with an 
appropriate generalization of homotopy, see in \cite{KW} and in 
\cite{D-V2}. It is given by the following lemma the proof of which is 
easy.

\begin{lemma}
Let $E$ be an $N$-differential module such that there are 
endomorphisms of
modules $h_k:E\rightarrow E$ for $k=0,1,\dots,N-1$ satisfying 
$\displaystyle{\sum^{N-1}_{k=0}}d^{N-1-k}h_kd^k=Id_E$; then one has
$H_{(n)}(E)=0$ for each integer $n$ with $1\leq n\leq N-1$.
\end{lemma}

In order to formulate the last criterion, we recall the definition of 
$q$-numbers. With $q\in\kb$, one associates a mapping $[.]_q:\mathbb 
N \rightarrow
\kb$, $n\mapsto [n]_q$,  which is defined by setting
$[0]_q=0$ and $[n]_q=1+\dots+q^{n-1}=\sum^{n-1}_{k=0} q^k$
for $n\geq 1,\ (q^0=1)$. For $n\in \mathbb N$ with $n\geq 1$, one 
defines the
$q$-factorial $[n]_q!\in \kb$ by $[n]_q\dots 1=\prod^n_{k=1}[k]_q$. 
For
integers $n$ and $m$ with $n\geq 1$ and $0\leq m\leq n$, one defines
inductively the $q$-binomial coefficients $\left[\begin{array}{c} n\\m
\end{array}\right]_q\in \kb$ by setting $\left[\begin{array}{c} n\\0
\end{array}\right]_q
= \left[\begin{array}{c} n\\n \end{array}\right]_q=1$ and
$\left[\begin{array}{c} n\\m \end{array}\right]_q +
q^{m+1}\left[\begin{array}{c} n\\m+1
\end{array}\right]_q=\left[\begin{array}{c} n+1\\m+1 
\end{array}\right]_q$ for
$0\leq m\leq n-1$. As in \cite{KW} let us introduce the following 
assumptions
$(A_0)$ and $(A_1)$ on the ring $\kb$ and the element $q$ of $\kb$ :\\
$(A_0)$\hspace{1cm} $[N]_q=0$\\
$(A_1)$\hspace{1cm} $[N]_q=0$ and $[n]_q$ is invertible for $1\leq 
n\leq N-1$,
$(n\in \mathbb N)$.\\
Notice that $[N]_q=0$ implies that $q^N=1$ and therefore that $q$ is
invertible. Furthermore if $q$ is invertible one has
$[n]_{q^{-1}}=q^{-n+1}[n]_q$, $\forall n\in \mathbb N$. Therefore 
Assumption
$(A_0)$, (resp. $(A_1)$), for $\kb$ and $q\in \kb$ is equivalent to 
Assumption
$(A_0)$, (resp. $(A_1)$), for $\kb$ and $q^{-1}\in \kb$. Let us give 
two typical examples:
\begin{enumerate}
\item
$\kb=\mathbb C,\ q\in \mathbb C$. Then Assumption $(A_0)$ means that 
$q$ is an
$N$-th root of  unity distinct of 1 and Assumption $(A_1)$ means that 
$q$ is
a primitive $N$-th root of  unity.
\item
$\kb=\mathbb Z_N=\mathbb Z/N\mathbb Z$, then $1\in \kb$ satisfies 
Assumption
$(A_0)$ and Assumption $(A_1)$ is satisfied if and only if $N$ is a 
prime
number.
\end{enumerate}

A useful result is that if $\kb$ and $q\in \kb$ satisfy Assumption 
$(A_1)$ then
one has $\left[\begin{array}{l}N\\m\end{array}\right]_q=0$ for $m\in
\{1,\dots,N-1\}$; notice that Assumption $(A_0)$ is not sufficient in 
order to
have this result.\\

We are now ready to state the last criterion \cite{D-V}.

\begin{lemma}
Suppose that $\kb$ and $q\in\kb$ satisfy $(A_1)$ and let $E$ be an
$N$-differential module. Assume that there is a module-endomorphism 
$h$ of $E$
such that $hd-qdh=Id_E$. Then one has $H_{(n)}(E)=0$ for each
integer $n$ with $1\leq n\leq N-1$.
\end{lemma}

In order to proof this lemma, one shows that in the unital 
$\kb$-algebra generated by $h$ and $d$ with the relation 
$hd-qdh=\bbbone$ one has 
$\displaystyle{\sum^{N-1}_{k=0}}d^{N-1-k}h^{N-1}d^k=[N-1]_q!\bbbone$, 
which implies the result in view of Lemma 3 since $[N-1]_q!$ is 
invertible in $\kb$ (see in \cite{KW} and in \cite{D-V2}).\\

It is obvious that Lemma 4 above is closely related to the theory of 
$q$-oscillators, (e.g. $d$ corresponds to the creation operator 
whereas $h$ corresponds to the annihilation operator), and this is 
the essence of the proof of \cite{D-V2}. As well known in physics, 
there is another natural way to produce creation operators with 
vanishing $N$-th powers which consists in considering parafermions of 
order $N-1$; this has the generalization we now describe.\\

As already pointed out (in Section 2 and Section 4) there is no 
natural tensor product between differential modules. The same is true 
for $N$-differential modules with $N$ fixed. However, if ($E',d'$) is 
an $N'$-differential module and if $(E'',d'')$ is an 
$N''$-differential module $(N',N''\geq 2)$ then one defines an 
$(N'+N''-1)$-differential $d$ on $E'\otimes E''$ by setting
\[
d=d'\otimes I'' + I'\otimes d''
\]
where $I'$ (resp. $I''$) denotes the identity mapping $Id_{E'}$ 
(resp. $Id_{E''}$) of $E'$ (resp. of $E''$). Therefore, a natural 
construction of an $N$-differential module consists in starting with 
$(N-1)$ ordinary differential modules $(E_i,d_i)$ and equipping their 
tensor product $E=E_1\otimes \cdots \otimes E_{N-1}$ with the 
$N$-differential 
\[
d=d_1\otimes I_2\otimes\cdots \otimes I_{N-1}+\cdots +I_1\otimes 
\cdots \otimes I_{N-2}\otimes d_{N-1}.
\]
If all the $(E_i,d_i)$ are identical, with $d_i$ being a fermionic 
creation operator, the above formula is the Green ansatz \cite{green} 
for the parafermionic creation operator of order $N-1$.\\

In the case where $\kb$ is a field, an $N$-differential module will 
be referred to as an $N$-{\sl differential vector space}. Assume that 
$E$ is a finite-dimensional $N$-differential vector space. Then one 
has $E\simeq \ker(d^n)\oplus \im (d^n)=Z_{(n)}(E)\oplus B_{(N-n)}(E)$ 
and $E\simeq \ker(d^{N-n})\oplus \im(d^{N-n})=Z_{(N-n)}(E)\oplus 
B_{(n)}(E)$ which together with $Z_{(n)}(E)\simeq B_{(n)}(E)\oplus 
H_{(n)}(E)$ and $Z_{(N-n)}(E)\simeq B_{(N-n)}(E)\oplus H_{(N-n)}(E)$ 
implies (since $\dim(E)<\infty$) that $H_{(n)}(E)$ {\sl and 
$H_{(N-n)}(E)$ are isomorphic}. In the case where $E$ is a 
finite-dimensional $N$-differential vector space over $\kb=\mathbb R$ 
or $\mathbb C$, one can show (see e.g. in \cite{Gre}) by decomposing 
$E$ into indecomposable factor for the action of the $N$-differential 
$d$ that one has an
isomorphism $E\simeq \oplus^N_{n=1}{\mathbf k}^n\otimes {\mathbf 
k}^{m_n}$,
$d\simeq \oplus^N_{n=2} D_n\otimes Id_{{\mathbf k}^{m_n}}$ with

\[
D_n=\left(
\begin{array}{ccccccc}
0 & 1 & 0 & . & . & . & 0\\
. & . & . & . &  & &.\\
. &   & . & . & . &  & .\\
. &   &   & . & . & . & .\\
. &  &   &   & .  & . & 0\\
. & & & & &.& 1\\
0 & . &. &. & . & . & 0\\
\end{array}
\right)\in M_n(\kb)
\]
where the {\sl multiplicities}\  $m_n$, $n\in \{1,\dots,N\}$, are 
invariants of
$(E,d)$ with\linebreak[4] $\sum^N_{n=1}n m_n=\dim (E)$. Notice that 
one has
$m_N\geq 1$ whenever $d^{N-1}\not=0$. Notice also that the above 
decomposition
of $d$ is its {\sl Jordan normal-form}. In terms of the 
multiplicities, one can
easily compute the dimensions of the vector spaces $H_{(k)}(E)$. The 
result is
given by the following proposition.

\begin{proposition}
Let $E$ be a finite dimensional $N$-differential vector space over 
$\mathbb R$
or $\mathbb C$ with multiplicities $m_n$, $n\in \{1,2,\dots,N\}$, 
then one has
for each integer $k$ with $1\leq k\leq N/2$
\[
\dim H_{(k)}(E)=\dim H_{(N-k)}(E)=\sum^k_{j=1}
\sum^{N-j}_{i=j}m_i.
\]
\end{proposition}
Although easy, that kind of results is useful for applications (see 
below).

\section{$N$-complexes}

An $N$-{\sl complex of modules}  \cite{Kap} or simply an $N$-{\sl 
complex} is an $N$-differential module $E$ which is $\mathbb 
Z$-graded, i.e. $E=\oplus_{n\in\mathbb Z}E^n$, with a homogeneous 
$N$-differential
$d$ of degree 1 or $-1$. When $d$ is of degree 1 then $E$ is referred 
to as {\sl a
cochain $N$-complex} and when $d$ is of degree $-1$ then $E$ is 
referred to as a
{\sl chain $N$-complex}. Here we adopt the cochain language and 
therefore in
the following an
$N$-{\sl complex}, without other specification, always means a cochain
$N\mbox{-complex}$ of modules. If $E$ is an $N\mbox{-complex}$ then 
the
$H_{(m)}(E)$ are $\mathbb Z$-graded
modules;
$H_{(m)}(E) = \oplusinf_{n\in \mathbb Z}H^n_{(m)}(E)$ with 
$H_{(m)}^n(E) = \ker
(d^m:E^n\rightarrow E^{n+m})/d^{N-m}(E^{n+m-N})$. In this case the 
hexagon
$(\calh^{\ell,m})$ of Lemma 1 splits into long exact sequences
$(\cals^{\ell,m}_p),\  p\in \mathbb Z$\\

\[
\begin{array}{lllll}
& & & &\dots \rightarrow H^{Nr+p}_{(m)}(E)\stackrel{[i]^\ell}{\hbox to
12mm{\rightarrowfill}} H^{Nr+p}_{(\ell+m)}(E) \stackrel{[d]^m}{\hbox 
to
12mm{\rightarrowfill}}H^{Nr+p+m}_{(\ell)}(E)\\
\\
(\cals^{\ell,m}_p) & && & \stackrel{[i]^{N-(\ell+m)}}{\hbox to
12mm{\rightarrowfill}}H^{Nr+p+m}_{(N-m)}(E)\stackrel{[d]^\ell}{\hbox 
to
12mm{\rightarrowfill}}H^{Nr+p+\ell+m}_{(N-(\ell+m))}(E)\\
\\
& && &\stackrel{[i]^m}{\hbox to
12mm{\rightarrowfill}}H^{Nr+p+\ell+m}_{(N-\ell)}(E)
\stackrel{[d]^{N-(\ell+m)}}{\hbox to
12mm{\rightarrowfill}}H^{N(r+1)+p}_{(m)}(E) \stackrel{[i]^\ell}{\hbox 
to
12mm{\rightarrowfill}} \dots
\end{array}
\]
\\
\noindent One has $(\cals^{\ell,m}_p)=(\cals^{\ell,m}_{p+N})$.\\

Let $E$ and $E'$ be $N$-complexes, a {\sl homomorphism of 
$N$-complexes} of
$E$ into $E'$ is a homomorphism of $N$-differential modules
$\varphi:E\rightarrow E'$ which is homogeneous of degree 0, (i.e.
$\varphi(E^n)\subset E^{\prime n}$). Such a homomorphism of 
$N$-complexes
induces module-homomorphisms $\varphi_\ast:H^n_{(m)}(E)\rightarrow
H^n_{(m)}(E')$ for $n\in \mathbb Z$ and $1\leq m\leq N-1$. Let 
$0\rightarrow E
\stackrel{\varphi}{\rightarrow} F \stackrel{\psi}{\rightarrow} 
G\rightarrow 0$
be a short exact sequence of $N$-complexes, then the hexagon 
$(\calh_n)$ of
Lemma 2 splits into  long exact sequences $(\cals_{n,p}), p\in\mathbb 
Z$\\

\[
\begin{array}{lllll}
& & & &\dots \rightarrow 
H^{Nr+p}_{(n)}(E)\stackrel{\varphi_\ast}{\hbox to
12mm{\rightarrowfill}} H^{Nr+p}_{(n)}(F) \stackrel{\psi_\ast}{\hbox to
12mm{\rightarrowfill}}H^{Nr+p}_{(n)}(G)\\
\\
(\cals_{n,p}) & && &
\stackrel{\partial}{\rightarrow}H^{Nr+p+n}_{(N-n)}(E)\stackrel{\varphi_\ast}{\hbox 
to 
12mm{\rightarrowfill}}H^{Nr+p+n}_{(N-n)}(F)\stackrel{\psi_\ast}{\hbox 
to 12mm{\rightarrowfill}} H^{Nr+p+n}_{(N-n)}(G)\\
\\
& && &\stackrel{\partial}{\rightarrow}H^{N(r+1)+p}_{(n)}(E)
\stackrel{\varphi_\ast}{\hbox to 12mm{\rightarrowfill}} \dots
\end{array}
\]
\\
\noindent One has again $(\cals_{n,p})=(\cals_{n,p+N})$.\\

In the following of this section, $(E^n)_{n\in\mathbb N}$ is a 
pre-cosimplicial module (see in Section 3), $E$ denotes the 
(positively) graded module $\oplusinf_n E^n$ and $q\in \kb$ is such 
that $[N]_q=0$, i.e. such that $\kb$ and $q\in \kb$ satisfy the 
assumption $(A_0)$ of Section 5.
One can construct a sequence $(d_n)_{n\in\mathbb N}$ of 
$N$-differentials of degree 1 on $E$ by using
$q\in\kb$ as above \cite{D-V2}. Here we
shall only consider the first two $d_0$ and $d_1$ which are the most 
natural ones. They are defined by
setting for $n\in \mathbb N$
\[
d_0=\sum^{n+1}_{i=0} q^i\fracf_i:E^n\rightarrow E^{n+1}
\]
and
\[
d_1=\sum^n_{i=0} q^i\fracf_i-q^n\fracf_{n+1}:E^n\rightarrow E^{n+1}.
\]
\begin{lemma}
One has $d^N_0=0$ and $d^N_1=0$.
\end{lemma}
This is a consequence of $[N]_q=0$ and of the relations $(\fracF)$; 
for a proof we
refer to \cite{D-V2}.\\
Thus $(E,d_0)$ and $(E,d_1)$ are $N$-complexes and, as shown in 
\cite{D-V2},
there are natural homomorphisms of the cohomology $H(E)$ of the
pre-cosimplicial module $(E^n)$ into the generalized cohomologies of 
these N-complexes.
In order to compute completely these generalized cohomologies we 
shall need
some more assumptions. We
shall need Assumption $(A_1)$ for $\kb$ and $q\in \kb$ and we shall 
restrict
attention to cosimplicial modules.
The generalized cohomologies of $(E,d_0)$ and $(E,d_1)$ are then 
given by the
following theorem \cite{D-V2}.

\begin{theo}
Let $\kb$ and $q\in\kb$ satisfy Assumption $(A_1)$ and let $(E^n)$ be 
a
cosimplicial
module. Then one has:
\[
\begin{array}{ll}
(0)& H^{Nr-1}_{(m)}(E,d_0)=H^{2r-1}(E),\ 
H^{N(r+1)-m-1}_{(m)}(E,d_0)=H^{2r}(E)\
\mbox{and}\\
& H^n_{(m)}(E,d_0)=0\ \ \mbox{otherwise},\\
\\
(1) & H^{Nr}_{(m)}(E,d_1)=H^{2r}(E), 
H^{N(r+1)-m}_{(m)}(E,d_1)=H^{2r+1}(E)\ \
\mbox{and}\\
& H^n_{(m)}(E,d_1)=0\ \ \mbox{otherwise},\\
\\
\end{array}
\]
for $r\in\mathbb N$ and $m\in\{1,\dots,N-1\}$.

\end{theo}
There is of course a dual statement for simplicial modules and the 
analogs of
$d_0$ and $d_1$ which are then of degree $-1$, see in \cite{D-V2}. 
The above
theorem and its dual version cover all the (co)simplicial cases 
investigated so
far that I know (\cite{May}, \cite{D-VK}, \cite{D-V}, \cite{KW} and
\cite{D-V2}). In \cite{D-V2} the generalized cohomology of $E$ for  
every  $d_n$ ($n\in \mathbb N$) was also computed in the case of a 
cosimplicial module $(E^n)$ as well as the generalized homologies of 
their chain analogs in the case of a simplicial module $(E_n)$ under 
assumption $(A_1)$ for $\kb$ and $q\in \kb$.
As a rule, we found
there that (in the (co)simplicial case) these generalized 
(co)homologies do
only depend on the ordinary (co)homology of the (co)simplicial 
module. In fact one of the ingredients in the proof of the above 
theorem is to use the whole sequence of $N$-differentials 
$(d_n)_{n\in\mathbb N}$ because, for any $p\in \mathbb N$ there is a 
$n_p\in\mathbb N$ such that $d_n$ coincides with the simplicial 
differential in degree $r$ (i.e. on $E^r$) whenever $n\geq n_p$ and 
$r\leq p$; the proof is nevertheless highly non trivial, (see in 
\cite{D-V2}).\\

Many notions for $N$-complexes do only depend on the underlying 
$\mathbb Z_N$-graduation ($\mathbb Z_N=\mathbb Z/N\mathbb Z$) so let 
us define a $\mathbb Z_N$-{\sl complex} to be an $N$-differential 
module which is $\mathbb Z_N$-graded with an $N$-differential which 
is homogeneous of degree 1. We have avoided the terminology $\mathbb 
Z_N$-$N$-complex since we shall not consider differential modules 
equipped with $\mathbb Z_N$ graduations with $N\geq 3$ and since a 
$\mathbb Z_2$-complex with the above definition is a $\mathbb 
Z_2$-complex according to the definition of Section 3. We now give an 
example of $\mathbb Z_N$-complex.\\

Let $\kb$ and $q\in \kb$ satisfy Assumption $(A_1)$ and let us 
introduce the
standard basis $E^k_\ell,\ (k,\ell\in\{1,\dots,N\}$), of the algebra 
$M_N(\kb)$
of $N\times N$ matrices defined by 
$(E^k_\ell)^i_j=\delta^k_j\delta^i_\ell$.
One has $E^k_\ell E^r_s=\delta^k_s E^r_\ell$ and $\sum^N_{n=1} 
E^n_n=\bbbone$.
It follows that one can equip $M_N(\kb)$ with a structure of $\mathbb
Z_N$-graded algebra, $M_N(\kb)=\oplus_{a\in\mathbb Z_N}M_N(\kb)^a$, 
by giving
to $E^k_\ell$ the degree $k-\ell\  \mbox{mod}(N)$. Let 
$e=\lambda_1E^2_1+\dots
+ \lambda_{N-1} E^N_{N-1}+\lambda_N E^1_N$ be an element of degree 1 
of
$M_N(\kb)$ and define the endomorphism $d$ by $d(A)=e A-q^a Ae$ for 
$A\in
M_N(\kb)^a$. One has $d^N=0$ so $(M_N(\kb),d)$ is a $\mathbb 
Z_N$-complex.
One verifies that $e^N=\lambda_1 \dots \lambda_N\bbbone$ and that 
$e^{N-1}
d(A)-qd(e^{N-1}A)=(1-q)\lambda_1\dots \lambda_N A$. Therefore if 
$1-q$ and the
$\lambda_i$ are invertible in $\kb$, Lemma 4 implies that
$H_{(n)}(M_N(\kb),d)=0$ for $n\in\{1,\dots,N-1\}$. It is worth 
noticing that
the above $N$-differential satisfies the {\sl graded $q$-Leibniz rule}
$d(AB)=d(A)B+q^a Ad(B)$, $\forall A\in M_N(\kb)^a$, $\forall B\in 
M_N(\kb)$.\\

It is clear that for any $N$-complex one has an underlying $\mathbb 
Z_N$-complex which is obtained by retaining only the degree modulo 
$N$. On the other hand starting from a $\mathbb Z_N$-complex 
$E=\oplusinf_{n\in\mathbb Z_N} E^n$ like $(M_N(\kb),d)$ above, one 
can construct an $N$-complex $\tilde E=\oplusinf_{n\in \mathbb 
Z}\tilde E^n$ by setting $\tilde E^n=E^{\pi(n)}$ where $\pi$ is the 
canonical projection of $\mathbb Z$ onto $\mathbb Z_N$, the 
definition of the $N$-differential on $\tilde E$ being obvious in 
terms of the one of $E$.\\

The content of Section 5 and Section 6 is based on \cite{D-V2} (see 
also in \cite{D-VK} and in \cite{D-V}). Particular $N$-complexes were 
introduced and analysed in \cite{May} for $\kb=\mathbb Z_N$, ($N$ 
prime). Several mathematicians wrote on $N$-complexes at the end of 
the 40's,  beginning of the 50's. The subject was reconsidered in 
\cite{Kap} and developed more recently in \cite{D-VK}, \cite{D-V}, 
\cite{KW} and \cite{D-V2}. In \cite{KW} an approach in the line of 
modern homological algebra \cite{cart:02} was developed with the 
introduction of generalizations of the functors Ext and Tor. 

\section{$N$-complexes of tensor fields}

In this section we shall describe $N$-complexes of tensor fields on 
$\mathbb R^D$ which generalize the complex $\Omega(\mathbb R^D)$ of 
differential forms \cite{D-VH}, \cite{D-VH2}. Therefore here the ring 
$\kb$ is the field $\mathbb R$ (or eventually $\mathbb C$ if one 
considers complex tensors). Furthermore, in such an $N$-complex, for 
each degree the tensor fields will be smooth mapping $x\mapsto T(x)$ 
of $\mathbb R^D$ into the vector space of covariant tensors of a 
given Young symmetry. Let us recall that this implies that the 
representation of $GL_D$ in the corresponding space of tensors is 
irreducible. For Young diagrams, etc. we refer to \cite{Ful} and for 
more details and developments we refer to \cite{D-VH}, \cite{D-VH2}.\\

Throughout the following $(x^\mu)=(x^1,\dots,x^D)$ denotes the
canonical coordinates of $\mathbb R^D$ and $\partial_\mu$ are the
corresponding partial derivatives which we identify with the
corresponding covariant derivatives associated to the canonical
flat linear connection of $\mathbb R^D$. Thus, for instance, if $T$
is a covariant tensor field of degree $p$ on $\mathbb R^D$ with
components $T_{\mu_1\dots\mu_p}(x)$, then $\partial T$ denotes
the covariant tensor field of degree $p+1$ with components
$\partial_{\mu_1}T_{\mu_2\dots\mu_{p+1}}(x)$. The operator
$\partial$ is a first-order differential operator which increases
by one the tensorial degree.\\

In this context, the space $\Omega(\mathbb R^D)$ of differential
forms on $\mathbb R^D$ is the graded vector space of (covariant)
antisymmetric tensor fields on $\mathbb R^D$ with graduation
induced by the tensorial degree whereas the exterior differential
$d$ is the composition of the above $\partial$ with
antisymmetrisation, i.e.
\[
d={\mathbf A}_{p+1}\circ \partial : \Omega^p(\mathbb R^D)\rightarrow
\Omega^{p+1}(\mathbb R^D)
\]
where ${\mathbf A}_p$ denotes the antisymmetrizer on tensors of 
degree $p$. One has  $d^2=0$ and
the Poincar\'e lemma asserts that the cohomology of the complex
$(\Omega(\mathbb R^D),d)$ is trivial, i.e. that one has
$H^p(\Omega(\mathbb R^D))=0$,  $\forall p\geq 1$ and 
$H^0(\Omega(\mathbb R^D))= \mathbb R$.\\

From the point of view of Young symmetry, antisymmetric
tensors correspond to Young diagrams (partitions) described 
by one column of cells, i.e. the space of values of
$p$-forms corresponds to one column of $p$ cells, $(1^p)$,
whereas ${\mathbf A}_p$ is the associated Young symmetrizer,
(see e.g. in \cite{Ful}).\\

There is a relatively easy way to generalize the pair
$(\Omega(\mathbb R^D),d)$ which we now describe. Let
$Y=(Y_p)_{p\in \mathbb N}$ be a sequence of Young diagrams
such that the number of cells of $Y_p$ is $p$, $\forall p\in
\mathbb N$ (i.e. such that $Y_p$ is a partition of the integer
$p$ for any $p$). We define $\Omega^p_Y(\mathbb R^D)$ to be
the vector space of smooth covariant tensor fields of degree
$p$ on $\mathbb R^D$ which have the Young symmetry type $Y_p$
and we let $\Omega_Y(\mathbb R^D)$ be the graded vector space
$\displaystyle{\oplusinf_p}\Omega^p_Y(\mathbb R^D)$. We then
generalize the exterior differential by setting $d={\mathbf
Y}\circ \partial$, i.e.
\[
d={\mathbf Y}_{p+1}\circ\partial :\Omega^p_Y(\mathbb
R^D)\rightarrow \Omega^{p+1}_Y(\mathbb R^D)
\]
where ${\mathbf Y}_p$ is now the Young symmetrizer on tensor
of degree $p$ associated to the Young symmetry $Y_p$. This $d$
is again a first order differential operator which is of
degree one, (i.e. it increases the tensorial degree by one),
but now, $d^2\not= 0$ in general. Instead, one has the
following result.

\begin{lemma}
Let $N$ be an integer with $N\geq 2$ and assume that $Y$ is
such that the number of columns of the Young diagram $Y_p$ is
strictly smaller than $N$ (i.e. $\leq N-1$) for any $p\in
\mathbb N$. Then one has $d^N=0$.
\end{lemma}

In fact the indices in one column are antisymmetrized and $d^N\omega$ 
involves necessarily at least two
partial derivatives $\partial$ in one of the columns since
there are $N$ partial derivatives involved and at most $N-1$
columns.\\

Thus if $Y$ satisfies the condition of Lemma 6,
$(\Omega_Y(\mathbb R^D),d)$ is an $N$-complex.
Notice that $\Omega^p_Y(\mathbb R^D)=0$ if the first column of
$Y_p$ contains more than $D$ cells and that therefore, if $Y$
satisfies the condition of Lemma 6, then $\Omega^p_Y(\mathbb
R^D)=0$ for $p>(N-1)D$. \\

One can also define a graded bilinear product on
$\Omega_Y(\mathbb R^D)$ by setting 
\[
(\alpha\beta)(x)={\mathbf Y}_{a+b}(\alpha(x)\otimes \beta(x))
\]
for $\alpha\in \Omega^a_Y(\mathbb R^D)$, $\beta\in
\Omega^b_Y(\mathbb R^D)$ and $x\in \mathbb R^D$. This product is
by construction bilinear with respect to the $C^\infty(\mathbb
R^D)$-module structure of $\Omega_Y(\mathbb R^D)$, 
($\Omega^0_Y(\mathbb R^D)=C^\infty(\mathbb R^D)$). However it is 
generically non associative.

In the following we shall not stay at this level of generality but,
for each $N\geq 2$ we shall choose a particular $Y$, denoted by
$Y^N=(Y^N_p)_{p\in\mathbb N}$, satisfying the condition of Lemma
6 which is maximal in the sense that all the rows are of maximal 
length $N-1$ except the last one (eventually). In other words the 
Young diagram with $p$ cells $Y^N_p$ is defined in the
following manner: write the division of $p$ by $N-1$, i.e. write
$p=(N-1)n_p+r_p$ where $n_p$ and $r_p$ are (the unique) integers
with $0\leq n_p$ and $0\leq r_p\leq N-2$ ($n_p$ is the quotient
whereas $r_p$ is the remainder), and let $Y^N_p$ be the Young
diagram with $n_p$ rows of $N-1$ cells and the last row with
$r_p$ cells (if $r_p\not= 0$). One has
$Y^N_p=((N-1)^{n_p},r_p)$, that is we fill the rows maximally.\\

\begin{figure}[h]
    \begin{center}
    \epsfxsize=4cm
\hspace{10cm}\epsfbox{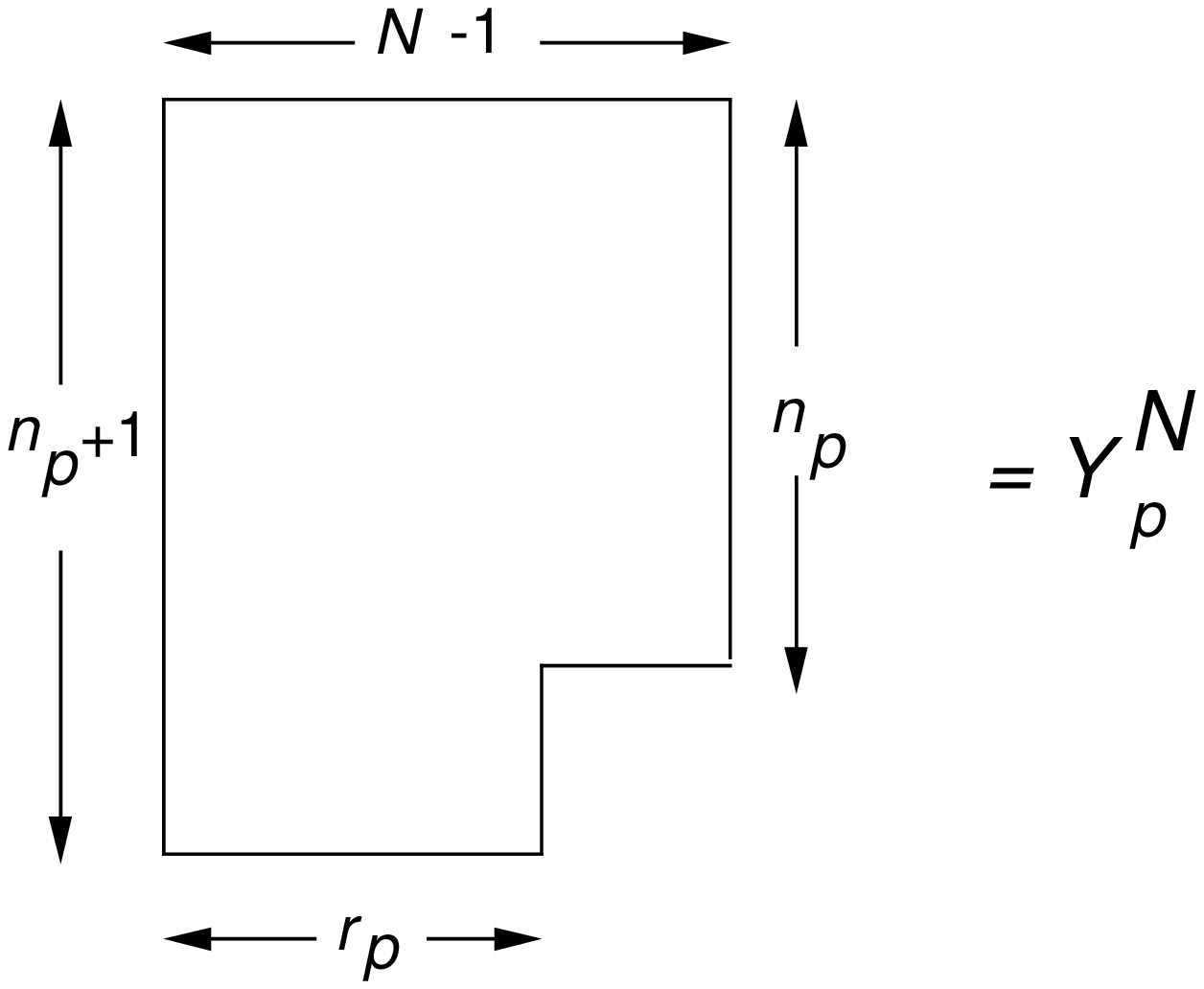}
    \end{center}
   \end{figure}

\noindent We shall denote $\Omega_{Y^N}(\mathbb R^D)$ and
$\Omega^p_{Y^N}(\mathbb R^D)$ by $\Omega_N(\mathbb R^D)$ and
$\Omega^p_N(\mathbb R^D)$. It is clear that $(\Omega_2(\mathbb
R^D),d)$ is the usual complex $(\Omega(\mathbb R^D),d)$ of
differential forms on $\mathbb R^D$. The $N$-complex
$(\Omega_N(\mathbb R^D),d)$ will be simply denoted by
$\Omega_N(\mathbb R^D)$. 
The Poincar\'e lemma admits the following generalization \cite{D-VH}, 
\cite{D-VH2}.

\begin{theo}
One has $H^{(N-1)n}_{(k)}(\Omega_N(\mathbb R^D))=0$, $\forall
n\geq 1$ and $H^0_{(k)}(\Omega_N(\mathbb R^D))$ is the space of
real polynomial functions on $\mathbb R^D$ of degree strictly
less than $k$ (i.e. $\leq k-1$) for $k\in\{1,\dots,N-1\}$.
\end{theo}

This statement reduces to the Poincar\'e lemma for $N=2$ but it
is a nontrivial generalization for $N\geq 3$ in the sense that,
the spaces $H^p_{(k)}(\Omega_N(\mathbb R^D))$
are nontrivial for $p\not=(N-1)n$ and in fact generically
infinite dimensional for $D\geq 3$, $p\geq N$.\\

The connection between the complex of differential forms on
$\mathbb R^D$ and the theory of classical gauge field of spin 1
is well known. Namely the subcomplex
\begin{equation}
\Omega^0(\mathbb R^D)\stackrel{d}{\rightarrow}\Omega^1(\mathbb
R^D)\stackrel{d}{\rightarrow}\Omega^2(\mathbb R^D)
\stackrel{d}{\rightarrow}\Omega^3(\mathbb R^D)
\label{eq1}
\end{equation}
has the following interpretation in terms of spin 1 gauge field
theory. The space $\Omega^0(\mathbb R^D)(=C^\infty(\mathbb R^D))$
is the space of infinitesimal gauge transformations, the space
$\Omega^1(\mathbb R^D)$ is 
the space of gauge potentials (which are the appropriate
description of spin 1 gauge fields to introduce local
interactions). The subspace $d\Omega^0(\mathbb R^D)$ of
$\Omega^1(\mathbb R^D)$ is the space of pure gauge configurations
(which are physically irrelevant), $d\Omega^1(\mathbb R^D)$ is
the space of field strengths or curvatures of gauge potentials.
The identity  $d^2=0$ ensures that the curvatures do not see the
irrelevant pure gauge potentials whereas, at this level, the
Poincar\'e lemma ensures that it is only these irrelevant
configurations which are forgotten when one passes from gauge
potentials to curvatures (by applying $d$). Finally $d^2=0$ also
ensures that curvatures of gauge potentials satisfy the Bianchi
identity, i.e. are in $\ker(d:\Omega^2(\mathbb R^D)\rightarrow
\Omega^3(\mathbb R^D))$, whereas at this level the Poincar\'e
lemma implies that conversely the Bianchi identity characterizes
the elements of $\Omega^2(\mathbb R^D)$ which are curvatures of
gauge potentials.\\

Classical spin 2 gauge field theory is the linearization of
Einstein geometric theory. In this case, and more generally in the 
linearization of (pseudo)riemannian geometry, the analog of
(\ref{eq1}) is a complex
$\cale^1\stackrel{d_1}{\rightarrow}\cale^2\stackrel{d_2}{\rightarrow}
\cale^3\stackrel{d_3}{\rightarrow}\cale^4$ where $\cale^1$ is the
space of covariant vector field $(x\mapsto X_\mu(x))$ on $\mathbb
R^D$, $\cale^2$ is the space of covariant symmetric tensor fields
of degree 2 ($x\mapsto h_{\mu\nu}(x))$ on $\mathbb R^D$,
$\cale^3$ is the space of covariant tensor fields of degree 4
$(x\mapsto R_{\lambda\mu,\rho\nu}(x))$ on $\mathbb R^D$ having the 
symmetries of
the Riemann curvature tensor  and where $\cale^4$
is the space of covariant tensor fields of degree 5 on $\mathbb
R^D$ having the symmetries of the left-hand side of the Bianchi
identity. The arrows $d_1, d_2, d_3$ are given by
\[
\begin{array}{l}
(d_1X)_{\mu\nu}(x)=\partial_\mu X_\nu(x)+\partial_\nu X_\mu (x)\\
(d_2h)_{\lambda\mu,\rho\nu}(x)=\partial_\lambda\partial_\rho 
h_{\mu\nu}(x)
+\partial_\mu\partial_\nu h_{\lambda\rho}(x)-\partial_\mu
\partial_\rho h_{\lambda\nu}(x) -\partial_\lambda\partial_\nu
h_{\mu\rho}(x)\\
(d_3R)_{\lambda\mu\nu,\alpha\beta}(x)=\partial_\lambda
R_{\mu\nu,\alpha\beta}(x)+\partial_\mu
R_{\nu\lambda,\alpha\beta}(x)+\partial_\nu
R_{\lambda\mu,\alpha\beta}(x).
\end{array}
\]
The symmetry of $x\mapsto R_{\lambda\mu,\rho\nu}(x)$,
$\left(\ \begin{tabular}{|c|c|}
\hline
$\lambda$ & $\rho$ \\
\hline
$\mu$ & $\nu$ \\
\hline
\end{tabular}\ 
\right)$, shows that $\cale^3=\Omega^4_3(\mathbb R^D)$ and that
$\cale^4=\Omega^5_3(\mathbb R^D)$; furthermore one canonically
has $\cale^1=\Omega^1_3(\mathbb R^D)$ and
$\cale^2=\Omega^2_3(\mathbb R^D)$. One also sees that $d_1$ and
$d_3$ are proportional to the 3-differential $d$ of
$\Omega_3(\mathbb R^D)$, i.e. $d_1\sim d:\Omega^1_3(\mathbb
R^D)\rightarrow \Omega^2_3(\mathbb R^D)$ and $d_3\sim d :
\Omega^4_3(\mathbb R^D)\rightarrow \Omega^5_3(\mathbb R^D)$.
The structure of $d_2$ looks different, it is of second order
and increases by 2 the tensorial degree. However it is easy to
see that it is proportional to $d^2:\Omega^2_3(\mathbb R^D)
\rightarrow \Omega^4_3(\mathbb R^D)$. Thus the analog of
(\ref{eq1}) is (for spin 2 gauge field theory)
\begin{equation}
\Omega^1_3(\mathbb R^D)\stackrel{d}{\rightarrow}
\Omega^2_3(\mathbb R^D)\stackrel{d^2}{\rightarrow}
\Omega^4_3(\mathbb R^D)\stackrel{d}{\rightarrow}
\Omega^5_3(\mathbb R^D)
\label{eq2}
\end{equation}
and the fact that it is a complex follows from $d^3=0$ whereas
the generalized Poincar\'e lemma (Theorem 3) implies that it is
in fact an exact sequence. Exactness at $\Omega^2_3(\mathbb
R^D)$ is $H^2_{(2)}(\Omega_3(\mathbb R^D))=0$ and exactness
at $\Omega^4_3(\mathbb R^D)$ is $H^4_{(1)}(\Omega_3(\mathbb
R^D))=0$, (the exactness at $\Omega^4_3(\mathbb R^D)$ is the
main statement of \cite{Gas}).\\

Thus what plays the role of the complex of differential forms
for the spin 1 (i.e. $\Omega_2(\mathbb R^D))$ is the 3-complex
$\Omega_3(\mathbb R^D)$ for the spin 2. More generally, for the
spin $S\in \mathbb N$, this role is played by the
$(S+1)$-complex $\Omega_{S+1}(\mathbb R^D)$. In particular, the
analog of Sequence (\ref{eq1}) for  spin 1 is the
complex
\begin{equation}
\Omega^{S-1}_{S+1}(\mathbb
R^D)\stackrel{d}{\rightarrow}\Omega^S_{S+1}(\mathbb
R^D)\stackrel{d^S}{\rightarrow}\Omega^{2S}_{S+1}(\mathbb
R^D)\stackrel{d}{\rightarrow} \Omega^{2S+1}_{S+1}(\mathbb R^D)
\label{eq3}
\end{equation}
for the spin $S$. The fact that (\ref{eq3}) is a complex was
known, \cite{dWF}, here it follows from $d^{S+1}=0$. One easily
recognizes that $d^S:\Omega^S_{S+1}(\mathbb R^D)\rightarrow
\Omega^{2S}_{S+1}(\mathbb R^D)$ is the generalized (linearized)
curvature of \cite{dWF}. Theorem 3 implies that sequence
(\ref{eq3}) is exact: exactness at $\Omega^S_{S+1}(\mathbb
R^D)$ is $H^S_{(S)}(\Omega_{S+1}(\mathbb R^D))=0$ whereas
exactness at $\Omega^{2S}_{S+1}(\mathbb R^D)$ is
$H^{2S}_{(1)}(\Omega_{S+1}(\mathbb R^D)=0$, (exactness at 
$\Omega^S_{S+1}(\mathbb R^D)$ was directly proved in \cite{DD} for 
the case $S=3$).\\

Finally, there is a generalization of Hodge duality for
$\Omega_N(\mathbb R^D)$, which is obtained by contractions of
the columns with the Kroneker tensor
$\varepsilon^{\mu_1\dots\mu_D}$ of $\mathbb R^D$ \cite{D-VH}, 
\cite{D-VH2}. When combined with Theorem 3, this
duality leads to another kind of results. A typical result of
this kind is the following one. Let $T^{\mu\nu}$ be a symmetric
contravariant tensor field of degree 2 on $\mathbb R^D$
satisfying $\partial_\mu T^{\mu\nu}=0$, (like e.g. the stress
energy tensor), then there is a contravariant tensor field
$R^{\lambda\mu\rho\nu}$ of degree 4 with the symmetry
$ \begin{tabular}{|c|c|}
\hline
$\lambda$ & $\rho$ \\
\hline
$\mu$ & $\nu$ \\
\hline
\end{tabular}\ 
$, (i.e. the symmetry of Riemann curvature tensor), such
that
\[
T^{\mu\nu}=\partial_\lambda\partial_\rho R^{\lambda\mu \rho\nu}
\]
In order to connect this result with Theorem 3, define
$\tau_{\mu_1\dots\mu_{D-1}
\nu_1\dots\nu_{D-1}}=\linebreak[4]T^{\mu\nu}\varepsilon_{\mu\mu_1\dots\mu_{D-1}}
\varepsilon_{\nu\nu_1\dots \nu_{D-1}}$. Then one has
$\tau\in\Omega^{2(D-1)}_3(\mathbb R^D)$ and conversely, any
$\tau\in \Omega^{2(D-1)}_3(\mathbb R^D)$ can be expressed in
this form in terms of a symmetric contravariant 2-tensor. It is
easy to verify that $d\tau=0$ (in $\Omega_3(\mathbb R^D))$ is
equivalent to $\partial_\mu T^{\mu\nu}=0$. On the other hand,
Theorem 3 implies that $H^{2(D-1)}_{(1)}(\Omega_3(\mathbb
R^D))=0$ and therefore $\partial_\mu T^{\mu\nu}=0$ implies that
there is a $\rho\in \Omega^{2(D-2)}_3(\mathbb R^D)$ such that
$\tau=d^2\rho$. The latter is equivalent to the above equation with
$R^{\mu_1\mu_2\ \nu_1\nu_2}$ proportional to
$\varepsilon^{\mu_1\mu_2\dots\mu_D}
\varepsilon^{\nu_1\nu_2\dots
\nu_D}\rho_{\mu_3\dots\mu_D\nu_3\dots\nu_D}$ and one verifies
that, so defined, $R$ has the correct symmetry. This result has been 
used in
\cite{W} in the investigation of the consistent deformations of
the free spin two gauge field action.

\section{Graded differential algebras and generalizations}

A {\sl graded differential algebra} is a (cochain) complex 
$\fraca=\oplus_{n\in \mathbb Z}\fraca^n$ with differential $d$ such 
that $\fraca$ is a $\mathbb Z$-graded associative unital 
$\kb$-algebra and such that $d$ is an antiderivation i.e. satisfies 
the graded Leibniz rule
\[
d(\alpha\beta)=d(\alpha)\beta+(-1)^a\alpha d(\beta)
\]
for any $\alpha\in\fraca^a$, $\beta\in\fraca$ and where 
$(\alpha,\beta)\mapsto \alpha\beta$ denotes the product of $\fraca$. 
If $\fraca$ is such a graded differential algebra with differential 
$d$, $\ker(d)$ is a graded unital subalgebra of $\fraca$ whereas 
$\im(d)$ is a graded two-sided ideal of $\ker(d)$ so the cohomology 
$H(\fraca)$ is a (unital associative) graded algebra. If $\fraca$ and 
$\fracB$ are two graded differential algebras, the tensor product 
$\fraca\otimes\fracB$ of the underlying complexes (as defined in 
Section 3) is again a graded differential algebra with product 
defined by
\[
(\alpha\otimes 
\beta)(\alpha'\otimes\beta')=(-1)^{ba'}\alpha\alpha'\otimes 
\beta\beta'
\]
for $\alpha\in \fraca$, $\beta\in\fracB^b$, $\alpha'\in\fraca^{a'}$ 
and $\beta'\in\fracB$. In the following, the product of a tensor 
product of graded algebras will be always the above one. With this 
convention, if $\kb$ is a field one has 
$H(\fraca\otimes\fracB)=H(\fraca)\otimes H(\fracB)$ for the 
corresponding cohomology algebras (which is the refined counterpart 
of Proposition 2 for graded differential algebras).\\

Let $(\fraca^n)_{n\in\mathbb N}$ be a pre-cosimplicial module (see in 
Section 3) such that $\fraca=\oplus \fraca^n$ is a (positively) 
graded algebra and assume that the cofaces homomorphisms $\fracf_i$ 
satisfy the following assumptions $(\fracM\fracF)$:

\noindent $(\fracM\fracF_1)\hspace{1.5cm} 
\fracf_i(\alpha\beta)=\left\{
\begin{array}{lll}
\fracf_i(\alpha)\beta &\mbox{if} & i\leq a\\
\alpha \fracf_{i-a}(\beta) &\mbox{if} & i>a
\end{array}
\right.
, i\in\{0,\dots,a+b+1\}$\\
and\\
$(\fracM\fracF_2)\hspace{1.5cm} \fracf_{a+1}(\alpha)\beta = \alpha 
\fracf_0(\beta)$\\
for $\alpha\in\fraca^a$ and $\beta\in\fraca^b$ (where 
$(\alpha,\beta)\mapsto \alpha\beta$ denote the product of $\fraca$). 
Then the corresponding complex $(\fraca,d)$ is a graded differential 
algebra. If furthermore $(\fraca^n)$ is a cosimplicial module with 
codegeneracy homomorphisms $\fracs_i$ satisfying the following 
assumption $(\fracM\fracS)$

\[
(\fracM\fracS)\hspace{1.5cm} \fracs_i(\alpha\beta)=\left\{
\begin{array}{lll}
\fracs_i(\alpha)\beta &\mbox{if} & i< a\\
\alpha \fracs_{i-a}(\beta) &\mbox{if} & i\geq a
\end{array}
\right.
\]
$i\in\{0,\dots,a+b-1\}$, then the subcomplex $N(\fraca)$ of 
normalized cochains of $\fraca$ is a graded differential subalgebra 
of $\fraca$. In \cite{D-V2}, a pre-cosimplicial module $(\fraca^n)$ 
as above with cofaces satisfying $(\fracM\fracF)$ (which was denoted 
there by $(\fraca\fracF)$) was called a {\sl pre-cosimplicial 
algebra} and in the case where $(\fraca^n)$ is furthermore a 
cosimplicial module with codegeneracies satisfying $(\fracM\fracS)$ 
(which was denoted there ($\fraca\fracS$)) it was called a {\sl 
cosimplicial algebra}, however it has been remarked by Max Karoubi 
that this terminology is misleading so we shall speak in the 
following of a $\fracM$-{\sl pre-cosimplicial module} in the first 
case and of a $\fracM$-{\sl cosimplicial module} in the second case. 
In fact $\fracM$-cosimplicial modules is what corresponds to graded 
differential algebras in an appropriate specific version of the 
Dold-Kan correspondence.\\

Let $\cala$ be an associative unital $\kb$-algebra and $\calm$ be a 
$(\cala,\cala)$-bimodule. As pointed out in Section 3, the 
$\calm$-valued Hochschild cochains give rise to a cosimplicial module 
$(C^n(\cala,\calm))_{n\in\mathbb N}$.
In the case $\calm=\cala$, $C(\cala,\cala)$ has a natural structure 
of $\mathbb N$-graded associative unital $\kb$-algebra with product 
$(\alpha,\beta)\mapsto
\alpha\beta$ given by
\[
\alpha\beta(x_1,\dots,x_{a+b})=\alpha(x_1,\dots,x_a)\beta(x_{a+1},\dots,x_{a+b}),
\]
for $\alpha\in C^a(\cala,\cala),\ \beta\in C^b(\cala,\cala), x_i\in 
\cala$.\\
It is easily verified that the assumptions $(\fracM\fracF)$ and 
$(\fracM\fracS)$ are satisfied so that $(C^n(\cala,\cala))$ is a 
$\fracM$-cosimplicial module. Thus $C(\cala,\cala)$ equipped with the 
simplicial (Hochschild) differential (as in Section 3) is a graded 
differential algebra and the submodule of normalized cochains is a 
graded differential subalgebra of $C(\cala,\cala)$.\\

Let again $\cala$ be an associative unital $\kb$-algebra and let us 
denote by
$\fracT(\cala)=\oplus_{n\in \mathbb N}\fracT^n(\cala)$ the tensor 
algebra over
$\cala$ of the $(\cala,\cala)$-bimodule $\cala\otimes \cala$. This is 
a (positively) graded associative
unital $\kb$-algebra with $\fracT^n(\cala)=\otimes^{n+1}\cala$ and 
product
$(x_0\otimes \dots\otimes x_n)(y_0\otimes \dots \otimes
y_m)=x_0\otimes\dots\otimes x_{n-1}\otimes x_ny_0\otimes y_1\otimes 
\dots
\otimes y_m$ for $x_i,y_j\in \cala$. One verifies that one defines a 
structure of $\fracM$-cosimplicial
module for $(\fracT^n(\cala))$ by setting\\
$\fracf_0(x_0\otimes\dots \otimes x_n)=\bbbone \otimes 
x_0\otimes\dots\otimes
x_n$\\
$\fracf_i(x_0\otimes\dots\otimes x_n)=x_0\otimes\dots \otimes 
x_{i-1}\otimes
\bbbone \otimes x_i\otimes\dots\otimes x_n\  \mbox{for}\ 1\leq i\leq 
n$\\
$\fracf_{n+1}(x_0\otimes \dots\otimes x_n)=x_0\otimes\dots\otimes 
x_n\otimes
\bbbone$\\
and\\
$\fracs_i(x_0\otimes\dots\otimes x_n)=x_0\otimes\dots\otimes 
x_ix_{i+1}\otimes
\dots\otimes x_n$ for $0\leq i\leq n-1$.\\
It follows that, equipped with the corresponding simplicial 
differential, $\fracT(\cala)$ is a graded differential algebra and 
that the submodule of normalized cochains is a graded differential 
subalgebra of $\fracT(\cala)$. This latter graded differential 
algebra will be denoted by $\Omega(\cala)$ and referred to as the 
{\sl universal graded differential envelope} of $\cala$ or simply the 
{\sl universal differential envelope} of $\cala$; it is characterized 
by the following universal property \cite{kar}, \cite{kar:02} (see 
also e.g. in \cite{RC} and \cite{D-V4}).

\begin{proposition} Any homomorphism $\varphi$ of unital
algebras of  $\cala$ into the subalgebra $\Omega^0$  of
elements of degree $0$ of a graded  differential algebra
$\Omega$  has a unique extension
$\tilde\varphi:\Omega(\cala)\rightarrow \Omega$ as a
homomorphism  of graded differential algebras.
\end{proposition}

The graded differential algebra $\Omega(\cala)$ is usually 
constructed in a different manner; the fact that it identifies with 
the graded differential algebra of normalized cochains of 
$\fracT(\cala)$ is well known. It is worth noticing here that 
$\Omega(\cala)$ is also the graded differential subalgebra of 
$\fracT(\cala)$ generated by $\cala$ (i.e. the smallest graded 
differential subalgebra which contains $\cala$).\\

We now come to an $N$-complex version of graded differential algebra 
$(N\geq 2)$. For that we shall need $q\in\kb$ such that Assumption 
$(A_1)$ of Section 5 is satisfied i.e. $[N]_q=0$ and $[n]_q$ 
invertible in $\kb$ for $n\in\{1,\cdots,N-1\}$. Throughout the 
following of this section, $N$ and $q\in\kb$ are fixed and such that 
$(A_1)$ is satisfied. The following lemma is basic for the 
generalization, \cite{D-V2}. In this lemma, (and in the following) 
$d_1$ is the $N$-differential defined in Section~6 for any 
pre-cosimplicial module.

\begin{lemma}
Suppose that $\kb$ and $q\in\kb$ satisfy Assumption $(A_1)$ and let 
$(\fraca^n)$ be a $\fracM$-pre-cosimplicial module. Then the 
$N$-differential $d_1$ satisfies the graded $q$-Leibniz rule, that is
\[
d_1(\alpha\beta)=d_1(\alpha)\beta+q^a\alpha d_1(\beta)
\]
for $\alpha\in\fraca^a$ and $\beta\in\fraca=\oplusinf_n\fraca^n$.
\end{lemma}

A unital associative graded algebra equipped with an $N$-differential 
satisfying the (above) {\sl graded $q$-Leibniz rule} will be referred 
to as a {\sl graded $q$-differential algebra} \cite{D-VK}, 
\cite{D-V2}. The content of the above lemma is that if $(\fraca^n)$ 
is a $\fracM$-pre-cosimplicial module then $(\fraca,d_1)$ is a graded 
$q$-differential algebra which is positively graded. If, furthermore 
$(\fraca^n)$ is a $\fracM$-cosimplicial module then the generalized 
cohomology of $(\fraca,d_1)$ is given in terms of the ordinary 
cohomology of $(\fraca^n)$ by Theorem 2.\\

Let $\cala$ be as above an associative unital algebra. It follows 
from Lemma 7 that $\fracT(\cala)$ equipped with the $N$-differential 
$d_1$ is a graded $q$-differential algebra (which is $\mathbb 
N$-graded). Let $\Omega_q(\cala)$ be the graded $q$-differential 
subalgebra of $\fracT(\cala)$ generated by $\cala$, i.e. the smallest 
subalgebra of $\fracT(\cala)$ which contains $\cala$ and which is 
stable by the $N$-differential $d_1$. As graded $q$-differential 
algebra, $\Omega_q(\cala)$ is characterized uniquely up to an 
isomorphism by the following universal property \cite{D-VK}, 
\cite{D-V2}.

\begin{proposition} 
Any homomorphism $\varphi$ of unital
algebras of  $\cala$ into the subalgebra $\Omega^0$  of
elements of degree $0$ of a graded  $q$-differential algebra
$\Omega$  has a unique extension
$\tilde\varphi:\Omega_q(\cala)\rightarrow \Omega$ as a
homomorphism  of graded $q$-differential algebras.
\end{proposition}

This is the $q$-analog of Proposition 5, (a homomorphism of graded 
$q$-differential algebra being a homomorphism of graded algebras 
permuting the $N$-differentials). For $N=2$, $\Omega_q(\cala)$ 
reduces to $\Omega(\cala)$. The graded $q$-differential algebra 
$\Omega_q(\cala)$ is referred to as the {\sl universal 
$q$-differential envelope} of $\cala$ \cite{D-VK}, \cite{D-V2}. The 
generalized cohomologies of $(\fracT(\cala),d_1)$ and of 
$\Omega_q(\cala)$ are generically trivial; one has the following 
result, \cite{D-V2}.

\begin{proposition}
Assume that $\cala$ admits a linear form $\omega\in \cala^\ast$ such 
that $\omega(\bbbone)=1$. Then the generalized cohomologies of 
$(\fracT(\cala),d_1)$ and of $\Omega_q(\cala)$ are given by:
\[
\begin{array}{llllll}
H^n_{(k)}(\fracT(\cala),d_1) & = & H^n_{(k)}(\Omega_q(\cala)) & = & 
0&  \mbox{for }\ \  n\geq 1\ \   \mbox{and}\\
H^0_{(k)}(\fracT(\cala),d_1) & = & H^0_{(k)}(\Omega_q(\cala)) & = & 
\kb,  & \forall k\in\{1,\cdots,N-1\}.
\end{array}
\]
\end{proposition}

Notice that the assumption of this proposition is satisfied if $\kb$ 
is a field and that the case $N=2$ means, under the same assumption, 
the triviality of the cohomologies of $\fracT(\cala)$ and 
$\Omega(\cala)$, (a well known fact, \cite{kar:02}).\\

The above discussion shows the naturality of the notion of graded 
$q$-differential algebra as ``$N$-generalization" or $q$-analog of 
the notion of graded differential algebra. This notion has a slight 
drawback which is the non existence of natural tensor products 
\cite{asi}; let us discuss this point.
It was shown in \cite{Kap} that if $q\in \kb$ is such
that Assumption $(A_1)$
is satisfied then one can construct a tensor product 
for\linebreak[4]  $N$-complexes in the
following manner. Let $(E',d')$ and $(E'',d'')$
be two $N$-complexes and let us define $d$ on $E'\otimes E''$ by 
setting
\[
d(\alpha'\otimes \alpha'')=d'(\alpha')\otimes
\alpha''+q^{a'}\alpha'\otimes d''(\alpha''), \ \forall \alpha'\in
E^{\prime a'},\ \forall \alpha''\in E'',
\]
one has by induction on $n\in \mathbb N$
\[
d^n(\alpha'\otimes \alpha'')=\sum^n_{m=0} q^{a'(n-m)}\left[ 
\begin{array}{l}
n\\ m \end{array}\right]_q d^{\prime m}(\alpha')\otimes d^{\prime 
\prime n-m}(\alpha''),
\]
therefore Assumption $(A_1)$ implies $d^N(\alpha'\otimes
\alpha'')=d^{\prime N}(\alpha')\otimes\alpha''+\alpha'\otimes 
d^{\prime\prime N}(\alpha'')=0$.
Unfortunately,
as pointed out in \cite{asi}, when $(E',d')$ and $(E'',d'')$ are 
furthermore
two graded $q$-differential algebras, $d$ fails to be a 
$q$-differential in
that it does not satisfy the graded $q$-Leibniz rule except for $q=1$ 
or
$q=-1$.\\

As for $N$-complexes, many notions for graded $q$-differential 
algebras do only depend on the underlying $\mathbb Z_N$-graduation so 
it is natural to consider the following $\mathbb Z_N$-graded version. 
A $\mathbb Z_N$-{\sl graded $q$-differential algebra} is a $\mathbb 
Z_N$-graded algebra equipped with a homogeneous endomorphism $d$ of 
degree 1 which is an $N$-differential, i.e. $d^N=0$, and which 
satisfies the graded $q$-Leibniz rule 
$d(\alpha\beta)=d(\alpha)\beta+q^a\alpha d(\beta)$ for $\alpha$ 
homogeneous of degree $a\in\mathbb Z_N$, (let us remind that $N$ and 
$q\in\kb$ are connected by Assumption $(A_1)$). We have already met 
such a $\mathbb Z_N$-graded $q$-differential algebra at the end of 
Section 6, (namely $M_N(\kb)$).\\

The notion of graded $q$-differential algebra was introduced in 
\cite{D-VK} for $\kb=\mathbb C$ as well as the construction of the 
universal $q$-differential envelopes. Here, we have followed the 
presentation of \cite{D-V2}.

\section{Subquotients and constraints}

Let $E$ be a module and let $u\in E^\ast$ be a linear form on $E$. To 
these data, one associates a graded differential algebra $K(u)$ which 
is constructed in the following manner. As an algebra, $K(u)$ is the 
exterior algebra (over $\kb$) $\wedge E$ of $E$ but it is equipped 
with the opposite graduation, i.e. $K(u)=\oplusinf_n K^n(u)$ with 
$K^n(u)=\wedge^{-n}E$ if $n\leq 0$ and $K^n(u)=0$ if $n>0$. The 
differential $d_u$ of $K(u)$ is then defined to be the unique 
homogeneous $\kb$-linear endomorphism of degree 1 of $K(u)$ 
satisfying the graded Leibniz rule and such that 
$d_u(e)=u(e)\in\kb=K^0(u)$ for any $e\in E=K^{-1}(u)$. One has 
$d^2_u=0$ so $K(u)$ is a graded differential algebra; the underlying 
complex is a {\sl Koszul complex} and will be referred to as the 
Koszul complex associated to the pair $(E,u)$.\\

Let $M$ be a smooth (finite-dimensional, connected, paracompact) 
manifold and let $V$ be a closed submanifold of $M$. The $\mathbb 
R$-algebra $C^\infty(M)$ of smooth functions on $M$ will play the 
role of $\kb$ and $I(V)$ will denote the ideal of a smooth function 
on $M$ vanishing on $V$. We introduce the following regularity 
assumption $(R_0)$ for the data $(M,V)$:\\
$(R_0)$ \begin{tabular}{l} 
 $I(V)$ is generated by $m$ functions $u_\alpha\in C^\infty(M)$, 
$\alpha\in\{1,\cdots,m\}$,  which\\
are independent on $V$ in the sense $du_1(x)\wedge\cdots\wedge 
du_m(x)\not= 0$, $\forall x\in V$.
\end{tabular}

\noindent Let $E=\mathbb R^m\otimesinf_{\mathbb R} C^\infty(M)$ be 
the free $C^\infty(M)$-module of rank $m$  with canonical basis 
denoted by $\pi_\alpha,\alpha\in \{1,\cdots,m\}$, and let $u\in 
E^\ast$ be the ($C^\infty(M)$-)linear form on $E$ defined by 
$u(\pi_\alpha)=u_\alpha\in C^\infty(M)$, for $\alpha\in 
\{1,\cdots,m\}$. The Koszul complex $K(u)$ associated to the pair 
$(E,u)$ identifies with the free graded-commutative unital 
$C^\infty(M)$-algebra generated by the $\pi_\alpha$ in degree $-1$ 
equipped with the differential $d_u$ above; This is a graded 
differential $C^\infty(M)$-algebra. Under Assumption $(R_0)$ one has 
with these notations the following result \cite{dv:2}.

\begin{lemma}
The cohomology $H(K(u))$ of $K(u)$ is given by $H^n(K(u))=0$ if 
$n\not= 0$ and $H^0(K(u))$ identifies canonically with the algebra 
$C^\infty(V)$ of smooth functions on $V$.
\end{lemma}

In fact, $d_u(K^{-1}(u))=I(V)$ so one has 
$H^0(K(u))=C^\infty(M)/I(V)$ which is canonically $C^\infty(V)$; 
notice that $C^\infty(V)$ is a $C^\infty(M)$-algebra. Notice also 
that, since the $\mathbb R$-algebra $C^\infty(M)$ is unital, $K(u)$ 
is also a graded differential algebra over $\mathbb R$.\\

Lemma 8 gives a homological description of the algebra of functions 
on a submanifold (under assumption $(R_0)$). Our aim is now to give a 
homological description of the algebra of functions on a quotient 
manifold and finally to mix both descriptions to obtain a homological 
description of the algebra of functions on a quotient of a 
submanifold (subquotient).\\

Let $V$ be a smooth manifold. Recall that a {\sl foliation} of $V$ is 
a vector subbundle $F$ of the tangent bundle $T(V)$  of $V$ which is 
such that the $C^\infty(V)$-module $\calf$ of sections of $F$ is also 
a Lie subalgebra of the Lie algebra of vector fields on $V$. In the 
following, we shall identify the foliation with $\calf$. The ideal 
$(\wedge\calf)^\perp$ of the algebra $\Omega(V)$ of differential 
forms on $V$ of forms vanishing on $\calf$ is a differential ideal so 
the quotient $\Omega(V)/(\wedge\calf)^\perp$ is a graded differential 
algebra over $\mathbb R$ which will be denoted by $\Omega(V,\calf)$ 
and referred to as the graded differential algebra of {\sl 
longitudinal forms} of $\calf$; its differential will be denoted by 
$d_\calf$. In fact $\Omega(V,\calf)$ is a subcomplex of the 
Chevalley-Eilenberg complex $C_\wedge(\calf,C^\infty(V))$, (see 
Appendix B), where $\calf$ acts by derivations on $C^\infty(V)$. Let 
$H(V,\calf)$ be the cohomology of longitudinal forms; in degree 0, 
$H^0(V,\calf)$ identifies with the $\mathbb R$-algebra of smooth 
functions on $V$ invariant by the action of the vector fields 
belonging to $\calf$. In the case where the quotient $V/\calf$ exists 
as smooth manifold and is such that the canonical projection 
$p:V\rightarrow V/\calf$ is a submersion,  $H^0(V,\calf)$ identifies 
with the algebra $C^\infty(V/\calf)$ of smooth functions on the 
quotient. Let us introduce the following regularity assumption 
$(R_1)$ for $(V,\calf)$:\\
$(R_1)$ \begin{tabular}{l}
$\calf$ is a free $C^\infty(V)$-module of rank $m'$, or equivalently\\
$F$ is a trivial vector bundle of rank $m'$.
\end{tabular}

\noindent If $(R_1)$ is satisfied $\Omega(V,\calf)$ identifies with 
the graded-commutative algebra\linebreak[4] 
$C^\infty(V)\otimesinf_{\mathbb R}\wedge \mathbb R^{m'}$; in fact, 
$\Omega(V,\calf)$ is then the free graded-commutative 
unital\linebreak[4] $C^\infty(V)$-algebra generated by the 
$\chi^{\alpha'}$ in degree 1, $\alpha'\in\{1,\cdots,m'\}$, where 
$(\chi^{\alpha'})$ is the dual basis of the  basis $(\xi_{\alpha'})$ 
of $\calf$. With these conventions, $d_\calf$ is given on the 
generators by
\begin{equation}
\left\{
\begin{array}{lll}
d_\calf f & = & \xi_{\alpha'}(f)\chi^{\alpha'},\ \  \forall f\in 
C^\infty(V)\\
d_\calf \chi^{\alpha'} & = & -\frac{1}{2} C^{\alpha'}_{\beta'\gamma'} 
\chi^{\beta'} \chi^{\gamma'}
\end{array}
\right.
\end{equation}
the $C^{\alpha'}_{\beta'\gamma'}\in C^\infty(V)$ being given by 
$[\xi_{\beta'},\xi_{\gamma'}]=C^{\alpha'}_{\beta'\gamma'} 
\xi_{\alpha'}$ i.e. 
$C^{\alpha'}_{\beta'\gamma'}=\chi^{\alpha'}([\xi_{\beta'},\xi_{\gamma'}])$. 
One must be aware of the fact that $\Omega(V,\calf)$ is a graded 
differential algebra over $\mathbb R$ (and not over $C^\infty(V)$).

It is worth noticing here that an infinite dimensional analog of the 
above appears in gauge theory; there, $V$ is replaced by the affine 
space of gauge potentials (connections), $\calf$ is replaced by the 
Lie algebra of the group of gauge transformations acting on gauge 
potentials whereas the analog of the $\chi^{\alpha'}$ are components 
of the ghost field $\chi(x)$. In this context the BRS differential 
\cite{BRS} corresponds to the longitudinal differential $d_\calf$, 
\cite{vial}, \cite{RS2} .\\

Let now $M$ be a smooth manifold, $V$ be a closed submanifold of $M$ 
and assume that $V$ is equipped with a foliation $\calf$. We want to 
combine the above constructions to produce a homological description 
of $C^\infty(V/\calf)$. More precisely, our aim is to produce a 
graded differential algebra which contains $C^\infty(M)$ and which 
has the longitudinal cohomology $H(V,\calf)$ as cohomology. We assume 
in the following that the assumption $(R_0)$ is satisfied by $(M,V)$ 
and that the assumption $(R_1)$ is satisfied by $(V,\calf)$. With 
$(M,V)$ satisfying $(R_0)$ is associated as above the Koszul complex 
$K(u)$ with differential $d_u$. Let 
$\calk=\oplusinf_{i,j}\calk^{i,j}$ be the bigraded algebra 
$K(u)\otimesinf_{\mathbb R} \wedge \mathbb R^{m'}$ with 
$\calk^{i,j}=K^i(u)\otimesinf_{\mathbb R}\wedge^j\mathbb R^{m'}$ i.e. 
$\calk^{i,j}=\wedge^{-i}\mathbb R^m\otimesinf_{\mathbb R} 
C^\infty(M)\otimesinf_{\mathbb R}\wedge^j\mathbb R^{m'}$ if $i\leq 
0\leq j$ and $\calk^{i,j}=0$ otherwise. One can also consider that 
$\calk$ is a $\mathbb Z$-graded algebra, $\calk=\oplusinf_n\calk^n$, 
for the total degree $\calk^n=\oplusinf_{i+j=n}\calk^{i,j}$. We shall 
again denote by $\pi_\alpha, \alpha\in\{1,\cdots, m\}$, and 
$\chi^{\alpha'}, \alpha'\in\{ 1,\cdots,m'\}$ the elements of $\calk$ 
corresponding to the canonical basis of $\mathbb R^m$ and of $\mathbb 
R^{m'}$. As graded $C^\infty(M)$-algebra, $\calk$ is the free 
graded-commutative unital $C^\infty(M)$-algebra generated by the 
$\pi_\alpha$ in degree $-1$ and the $\chi^{\alpha'}$ in degree 1. One 
recovers the bidegree by giving the bidegree $(-1,0)$ to the 
$\pi_\alpha$ and the bidegree $(0,1)$ to the $\chi^{\alpha'}$. Let us 
extend the differential $d_u$ of $K(u)$ as the unique antiderivation 
$\delta_0$ of $\calk$ such that $\delta_0\chi^{\alpha'}=0$, $\delta_0 
f=0$ for $f\in C^\infty(M)$ and $\delta_0\pi_\alpha=u_\alpha$; one 
still has $\delta^2_0=0$ so $\calk$ equipped with $\delta_0$ is a 
graded differential algebra. Furthermore since $\delta_0$ is 
homogeneous for the bidegree (of bidegree $(1,0)$) the cohomology 
$H(\delta_0)$ of $(\calk,\delta_0)$ is bigraded, 
$H(\delta_0)=\oplusinf_{i,j}H^{i,j}(\delta_0)$, and Lemma 8 implies 
that $H^{i,j}(\delta_0)=0$ if $i\not= 0$ and that 
$H^{0,j}(\delta_0)=C^\infty(V)\otimesinf_{\mathbb R} \wedge^j\mathbb 
R^{m'}$ in other words one has the following lemma.

\begin{lemma}
As a graded $\mathbb R$-algebra, the cohomology $H(\delta_0)$ of 
$(\calk,\delta_0)$ identifies with the graded algebra 
$\Omega(V,\calf)$ of longitudinal forms.
\end{lemma}

The following lemma states that there is an antiderivation of $\calk$ 
which induces the longitudinal differential $d_\calf$ on 
$H(\delta_0)$.

\begin{lemma}
There is an antiderivation of $\delta_1$ of degree 1 of $\calk$ which 
is homogeneous for the bidegree of bidegree $(0,1)$, which satisfies 
$\delta_0\delta_1+\delta_1\delta_0=0$ and which induces the 
longitudinal differential $d_\calf$ on $H(\delta_0)=\Omega(V,\calf)$.
\end{lemma}

\noindent\underbar{Proof}. The longitudinal differential is given by 
(4) on $C(V)\otimesinf_{\mathbb R} \wedge \mathbb R^{m'}$. It follows 
from our assumptions that there are vector fields $\tilde 
\xi_{\alpha'}$ on $M$ such that their restrictions to $V$ are tangent 
to $V$ and coincide with the $\xi_{\alpha'}$, 
$\tilde\xi_{\alpha'}\restriction V=\xi_{\alpha'}$ for 
$\alpha'\in\{1,\cdots,m'\}$. Similarily there are $\tilde 
C^{\alpha'}_{\beta'\gamma'}\in C^\infty(M)$ such that $\tilde 
C^{\alpha'}_{\beta'\gamma'}\restriction 
V=C^{\alpha'}_{\beta'\gamma'}\in C^\infty(V)$. Define then $\delta_1$ 
on $C^\infty(M)\otimes_{\mathbb R}\wedge \mathbb R^{m'}$ by

\begin{equation}
\left\{
\begin{array}{lll}
\delta_1 f & = & \tilde \xi_{\alpha'}(f)\chi^{\alpha'},\ \ \forall 
f\in C^\infty(M)\\
\delta_1\chi^{\alpha'} & = & -\frac{1}{2} \tilde 
C^{\alpha'}_{\beta'\gamma'}\chi^{\beta'}\chi^{\gamma'}
\end{array}
\right.
\end{equation}

One has $(\delta_0\delta_1+\delta_1\delta_0)f=\delta_0\delta_1 f=0$ 
and 
$(\delta_0\delta_1+\delta_1\delta_0)\chi^{\alpha'}=\delta_0\delta_1\chi^{\alpha'}=0$ 
for $f\in C^\infty(M)$ and $\alpha'\in\{ 1,\cdots,m'\}$. On the other 
hand, one has 
$\delta_1\delta_0\pi_\alpha=\delta_1u_\alpha=\tilde\xi_{\alpha'}(u_\alpha)\chi^{\alpha'}$ 
and, by construction $\tilde\xi_{\alpha'}(u_\alpha)$ vanishes on $V$ 
so $\tilde\xi_{\alpha'}(u_\alpha)=A_{\alpha'\alpha}^\beta u_\beta$ 
for some $A^\beta_{\alpha'\alpha}\in C^\infty(M)$. By setting 
$\delta_1\pi_\alpha=-A^\beta_{\alpha'\alpha} \pi_\beta\chi^{\alpha'}$ 
and by extending $\delta_1$ to $\calk$ by the antiderivation 
property, one has $\delta_0\delta_1+\delta_1\delta_0=0$ so $\delta_1$ 
induces an antiderivation of degree 1 of 
$H(\delta_0)=C(V)\otimesinf_{\mathbb R}\wedge \mathbb R^{m'}$ which 
coincides with $d_\calf$ in view of (4). $\square$\\
As it is apparent in Formula (5), $\delta_1$ is an antiderivation of 
$\calk$ considered as a graded algebra over $\mathbb R$ (and not over 
$C^\infty(M)$ in contrast with $\delta_0$).

\begin{lemma}
There are antiderivations $\delta_r$ of degree 1 of the graded 
$\mathbb R$-algebra $\calk$ with $\delta_r$ homogeneous for the 
bidegree of bidegree $(1-r,r)$ for $r\geq 2$, such that one has with 
$\delta_0$ and $\delta_1$ as above $\sum_{r+s=n}\delta_r\delta_s=0$ 
for any $n\in \mathbb N$.
\end{lemma}

For the proof we refer to the proof of Theorem 3.7 of \cite{dv:2}. 
This is a proof by induction on $n$ using $H^{1-r,r+1}(\delta_0)=0$ 
and $H^{1-r,r+2}(\delta_0)=0$ for $r\geq 2$. Notice that $\delta_r=0$ 
if $r>m'$ or $r>m+1$.

\begin{theo}
Let $\delta_r$ $(r\geq 0)$ be as above then 
$\delta=\sum_{r\geq0}\delta_r$ is a differential of the graded 
$\mathbb R$-algebra $\calk$ and the cohomology $H(\delta)$ of the 
graded differential algebra $(\calk,\delta)$ identifies with the 
longitudinal cohomology $H(V,\calf)$.
\end{theo}

Again we refer to \cite{dv:2} (the proof of Theorem 3.8 there); the 
first part of the statement is obvious, the identification of 
$H(\delta)$ with $H(V,\calf)$ follows essentially from an elementary 
spectral sequence argument.\\

Let $(M,\omega)$ be a {\sl symplectic manifold} (i.e. a smooth 
manifold $M$ equipped with a closed nondegenerate 2-form $\omega$) 
and let $V$ be a closed submanifold of $M$. We denote by $\omega_V$ 
the closed 2-form $i^\ast(\omega)$ on $V$ induced by the inclusion 
$i:V\rightarrow M$. In general $\omega_V$ is degenerate; its 
{characteristic distribution} $F$ is the set of tangent vectors $X$ 
of $V$ such that $i_X\omega=0$. It follows from the equation 
$d\omega_V=0$ that the $C^\infty(V)$-module $\calf$ of vector fields 
on $V$ which are valued in $F$ is a Lie subalgebra of the Lie algebra 
of vector fields. Therefore if $\omega_V$ is of constant rank, which 
will be assumed in the sequel, $\calf$ is a foliation of $V$. In fact 
we shall assume not only that $\omega_V$ is of constant rank but also 
that the quotient $V/\calf=M_0$ is a smooth manifold and that the 
canonical projection $p:V\rightarrow M_0$ is a submersion. With these 
regularity assumptions, $\omega_V$ has a projection $\omega_0$ on 
$M_0$ which is, by construction, a closed nondegenerate 2-form. Thus 
$(M_0,\omega_0)$ is a symplectic manifold which is referred to as 
{\sl the reduced phase space} and which is the natural phase space 
for a hamiltonian system on $M$ which is constrained to move on $V$. 
One has $i^\ast(\omega)=\omega_V=p^\ast(\omega_0)$. The algebra of 
observables of such a constrained system is $C^\infty(M_0)$ which 
identifies with the longitudinal cohomology of degree 0, 
$C^\infty(M_0)=H^0(V,\calf)$. Thus if $(M,\omega)$ and $V$ are such 
that $(R_0)$ is satisfied for $(M,V)$ and $(R_1)$ is satisfied for 
$(V,\calf)$, one can use Theorem 4 to compute $C^\infty(M_0)$ and 
more generally $H(V,\calf)$. The graded differential algebra $\calk$ 
is the ghost complex appropriate to the situation and $\delta$ is the 
corresponding BRS differential.\\

The specificity of the above situation is that $\calf$ does only 
depend on the submanifold $V$ of the symplectic manifold 
$(M,\omega)$; in particular if assumptions $(R_0)$ and $(R_1)$ are 
satisfied, one can show easily that $m\geq m'$ and, on the other hand 
$m+m'=\dim (M)-\dim(M_0)$ is necessarily even since $M$ and $M_0$ are 
both symplectic (and finite-dimensional). The case where the ideal 
$I(V)$ of smooth functions on $M$ which vanish on $V$ is stable by 
the Poisson bracket (associated to $\omega)$ is referred to as the 
case of {\sl first class constraints} or the {\sl coisotropic case}. 
In such a case, one has $m=m'$ and Assumption $(R_0)$ implies 
$(R_1)$; indeed in this case with Assumption $(R_0)$ the hamiltonian 
vector fields $\ham(u_\alpha)$ of the $u_\alpha$ have restrictions to 
$V$ which are tangent to $V$ and form a basis of the 
$C^\infty(V)$-module $\calf$, (see e.g. in \cite{dv:2}). This case 
has the further property that one can extend the Poisson bracket in a 
superbracket on $\calk$ by setting 
$\{\pi_\alpha,\chi^\beta\}=\delta^\beta_\alpha$,  
$\{\chi^\alpha,\chi^\beta\}=\{\pi_\alpha,\pi_\beta\}=0$, 
$\{\pi_\alpha,f\}=\{\chi^\alpha,f\}=0$ for $f\in C^\infty(M)$ and 
that $\calk$ can then be interpreted as the algebra of ``functions" 
on a ``super phase space". Moreover in this case the BRS differential 
$\delta$ can be realized as superhamiltonian, i.e. $\delta\varphi=\{ 
Q,\varphi\}$, $\forall \varphi\in\calk$, for some $Q\in \calk$ of 
total degree 1, \cite{he:1}.  In this case it has been shown in 
\cite{he:1} that the arbitrariness of the whole construction is a 
canonical transformation of the super phase space.\\

In gauge theory, the usual ghost complex without antighosts was 
understood early as a Lie algebra cochain complex (see e.g. in 
\cite{bc-r}, \cite{RS}, \cite{dvtv}, \cite{dv:1}) or as a complex of 
longitudinal forms (see in \cite{vial}, \cite{RS2}). This led through 
the Koszul formula \cite{kosz:00} to the interpretation of the 
corresponding ghosts as components of the Maurer-Cartan form of the 
gauge group (see also in \cite{sul}).
 The key of the understanding of the antighost or conjugate ghost 
part in terms of usual mathematical concepts appears in \cite{McM} 
where it was shown that they provide Koszul resolutions. This led to 
the homological approach to constrained systems developed e.g. in 
\cite{ros}, \cite{sta}, \cite{dv:2} and \cite{KS} in terms of 
standard mathematical objects which is partly described above.\\

Assumption $(R_0)$ means regular submanifold $V$. One can generalize 
the above constructions in several directions without such a 
regularity. In the case of the first class constrained hamiltonian 
systems this has been investigated in \cite{sta} and in \cite{fhst} 
where BRS cohomology with ghosts of ghosts has been applied.
Finally, it is worth noticing here that an ``infinite dimensional" 
form of Theorem 4 applies directly to the antifield formalism 
\cite{Fhe}, \cite{HT} and is also implicit behind the ghost 
lagrangian formalism of gauge theory \cite{FP}, \cite{BRS}.
 
\section{$N$-complex versions of BRS methods}

The canonical approach to the quantum 
Wess-Zumino-Novikov-Witten\linebreak[4] (WZNW) model gives rise to a 
finite-dimensional quantum group gauge problem for the zero modes. 
This has been studied in a convenient form for us in \cite{fht1}, 
\cite{fht2}. The result is a finite-dimensional gauge model in which 
the physical state space appears as a quotient $\calh'/\calh''$ where 
$\calh'$ is a subspace of the original finite-dimensional indefinite  
metric space whereas $\calh''$ is the subspace of ``null vectors" 
(isotropic subspace) of $\calh'$. Using the results of \cite{fht1}, 
\cite{fht2}, it was shown in \cite{D-VT1} that, in the case of the 
$SU(2)$ WZNW model, the physical state space can be realized as a 
direct sum $\oplusinf_{n=1}^{N-1}H_{(n)}(\calh_I,A)$ where 
$(H_{(n)}(\calh_I,A))$ is the generalized homology of an 
$N$-differential vector space $\calh_I$ with $N$-differential $A$. In 
fact, for the level $k$ representation of the 
$\widehat{\fracs\fracu}(2)$ Kac-Moody algebra, $A$ satisfies $A^N=0$  
with $N=k+2$. The $N^4$-dimensional space $\calf\otimes \bar 
\calf=\calh$ of chiral zero modes carries a representation of the 
quantum group $U_q(\fracs\fracl_2)\otimes U_q(\fracs\fracl_2)$ where 
$q=e^{i(\pi/N)}$; it is a representation of the usual 
finite-dimensional quotient $\calu_q$ of $U_q(\fracs\fracl_2)\otimes 
U_q(\fracs\fracl_2)$ at the primitive root of unity $q$ ($q^{2N}=1$). 
The $N$-differential $A$ of $\calh$ commutes with the action of the 
Hopf algebra $\calu_q$ so the $(2N-1)$-dimensional subspace $\calh_I$ 
of $\calu_q$-invariant vectors is stable by $A$ and it is the 
generalized homology of the $N$-differential vector space 
($\calh_I,A$) which is of interest. In \cite{D-VT2} we produced an 
$N$-differential vector space which contains $\calh$ and has the same 
generalized homology as ($\calh_I,A$). It is this construction which 
will be explained in a very general setting in what follows.\\

In short, one has a vector space $\calh$ on which act a Hopf
algebra $\calu_q$ and a nilpotent endomorphism $A$
satisfying $A^N=0$. The action of the algebra $\calu_q$ commutes
with $A$, i.e. one has on $\calh$ : $ [A,X]=0,\ \ \ \forall 
X\in \calu_q$.
It follows that the subspace $\calh_I$ of $\calu_q$-invariant
vectors in $\calh$ is stable by $A$, i.e. $A(\calh_I)\subset
\calh_I$. Thus $(\calh_I,A)$ is an $N$-differential subspace of
the $N$-differential vector space $(\calh,A)$ and it turns out
that the ``interesting object" (the physical space) is the
generalized homology of $(\calh_I,A)$. We would like to avoid the
restriction to the invariant subspace $\calh_I$ that is, in
complete analogy with the BRS methods, we would like to define an
extended $N$-differential space in such a way that the
$\calu_q$-invariance is captured by its $N$-differential in the
sense that it
has the same generalized homology as $(\calh_I,A)$.\\
The most natural thing to do is to try to construct a nilpotent
endomorphism $Q$ of $\calh$ with $Q^N=0$ such that its
generalized homology coincides with the one of $A$ on $\calh_I$
i.e. such that one has $H_{(n)}(\calh,Q)=H_{(n)}(\calh_I,A),\ \ 
\forall n\in \{1,\dots,N-1\}$. It turns out that this is impossible 
in general. Indeed in the above case (for the $SU(2)$ WZNW model)  
$\calh$ is finite dimensional and then Proposition 4 (see in Section 
5) imposes strong constraints connecting $\dim\calh$ and the $\dim 
H_{(n)}(\calh,Q)=\dim H_{(n)}(\calh_I,A)$ for $n\in\{1,\cdots,N-1\}$ 
which  are not satisfied \cite{D-VT2}. This is not astonishing since 
in the usual BRS methods one has to add the ghost sector (see e.g. in 
last section or in Section 4).\\

We first present an abstract optimal construction in which the Hopf 
algebra $\calu_q$ plays no role. We assume that $(\calh,A)$ is an 
$N$-differential vector space, that there is a subspace $\calh_I$ of 
$\calh$ stable by $A$ and we shall construct an $N$-differential 
vector space $(\calh^\bullet,Q)$ with $\calh\subset \calh^\bullet$ 
such that $H_{(n)}(\calh^\bullet,Q)=H_{(n)}(\calh_I,A)$ for all $n\in 
\{ 1,\cdots,N-1\}$. Throughout the following $q^2$ is still a 
primitive $N$-th root of  unity $(q^N=-1)$. Let us define the graded 
vector space
$\calh^\bullet=\displaystyle{\oplusinf_{n\geq 0}}\calh^n$ by
$\calh^0=\calh$, $\calh^n=\calh/\calh_I$ for $1\leq n\leq N-1$
and $\calh^n=0$ for $n\geq N$. One then defines an endomorphism
$d$ of degree 1 of $\calh^\bullet$ by setting
$d=\pi:\calh^0\rightarrow \calh^1$ where $\pi:\calh\rightarrow
\calh/\calh_I$ is the canonical projection,
$d=\id:\calh^n\rightarrow \calh^{n+1}$ for $1\leq n\leq N-2$
where $\id$ is the identity mapping of $\calh/\calh_I$ onto
itself and $d=0$ on $\calh^n$ for $n\geq N-1$. One has $d^N=0$
and therefore $(\calh^\bullet,d)$ is an $N$-complex, so its 
generalized (co)homology is
graded  $H_{(k)}(\calh^\bullet,d)=
\displaystyle{\oplusinf_{n\geq 0}} H^n_{(k)}(\calh^\bullet,d)$. 
It is given by the following easy lemma.
\begin{lemma}
One has $H^n_{(k)}(\calh^\bullet,d)=0$ for $n\geq 1$
and
$H^0_{(k)}(\calh^\bullet,d)=\calh_I$, $\forall k\in \{1,\dots,
N-1\}$.
\end{lemma}

It is worth noticing here that given the vector space $\calh$
together with the subspace $\calh_I$, the $N$-complex
$(\calh^\bullet,d)$ is characterized (uniquely up to an
isomorphism) by the following universal property (the proof of
which is straightforward).
\begin{lemma}
Any linear mapping $\alpha:\calh\rightarrow \calc^0$ of $\calh$
into the subspace $\calc^0$ of elements of degree $0$ of an
$N$-complex $(\calc^\bullet,d)$ which satisfies $d\circ
\alpha(\calh_I)=0$ extends uniquely as a
homomorphism
$\bar\alpha:(\calh^\bullet,d)\rightarrow (\calc^\bullet,d)$ of
$N$-complexes.
\end{lemma}
By using this universal property one can extend $A$ to 
$\calh^\bullet$ in the following manner.

\begin{lemma}
The endomorphism $A$ of $\calh=\calh^0$ has a unique extension
to $\calh^\bullet$, again denoted by $A$, as a homogeneous
endomorphism of degree $0$ satisfying $Ad-q^2\ dA=0$. On
$\calh^\bullet$, one has $A^N=0$ and $(d+A)^N=0$.
\end{lemma}

Thus $Q=d+A$ is an $N$-differential on $\calh^\bullet$ and we have 
the following result.

\begin{theo}
The generalized $Q$-homology of $\calh^\bullet$ coincides with
the generalized $A$-homology of $\calh_I$, i.e. one has 
$H_{(k)}(\calh^\bullet,Q)=H_{(k)}(\calh_I,A)$
for $1\leq k\leq N-1$.
\end{theo}

Notice that $(\calh^\bullet,Q)$ is only an $N$-differential vector 
space  {\sl and not} an $N$-complex since $d+A=Q$ is inhomogeneous.\\

In the problem of the zero modes of the $SU(2)$ WZNW model, $\calh_I$ 
is the invariant subspace of $\calh$ by the action of the quantum 
group (i.e. the Hopf algebra) $\calu_q$ which plays the role of a 
gauge group, or more precisely of the universal enveloping algebra of 
the Lie algebra of a gauge group, so (in view of Theorem 1) it is 
natural to produce a construction where $\calu_q$ and its 
(Hochschild) cohomology enter as in the usual BRS construction for 
gauge theory in order to get a similar ``geometrico-physical" 
interpretation. This is the aim of the end of this section. Since the 
above construction based on universal property is quite minimal, one 
cannot be astonished that it occurs as an $N$-differential subspace 
of the following one.

By definition $\calh_I$ is the set of $\Psi\in \calh$ such that 
$X\Psi=\Psi\varepsilon(X)$  for any\linebreak[4] $X\in\calu_q$, where 
$\varepsilon$ denotes the counit of $\calu_q$. This means that if one 
considers $\calh$ as a\linebreak[4] $(\calu_q,\calu_q)$-bimodule by 
equipping it with the trivial right action given by the counit, 
$\calh_I$ identifies with the $\calh$-valued Hochschild cohomology in 
degree 0 of $\calu_q$, i.e. $\calh_I=H^0(\calu_q,\calh)$. The idea of 
the construction is to mix the Hochschild differential with $A$ in a 
similar way as the mixing of $\delta_0$ with $d_\calf$ in last 
section. However, $A$ is an $N$-differential whereas the Hochschild 
differential is an ordinary  differential i.e. a 2-differential. 
Fortunately the next lemma shows that for the description of 
$\calh_I$ one can replace the Hochschild differential by the 
$N$-differential $d_1$ of Section 6 with the replacement of $q$ by 
$q^2$ since here it is $q^2$ which is a primitive $N$-th root of 
unity. To simplify the notations, this $N$-differential $d_1$ on 
$C(\calu_q,\calh)$ will be denoted by $d$. That is the 
$N$-differential $d$ is defined by 
\[
\begin{array}{lll}
d(\omega)(X_0,\dots,X_n) & = & X_0\omega(X_1,\dots,X_n)\nonumber\\
& + & \sum^n_{k=1}
q^{2k}\omega(X_0,\dots,(X_{k-1}X_k),\dots,X_n)\nonumber\\
& - & q^{2n}\omega(X_0,
\dots, X_{n-1})\varepsilon(X_n).
\end{array}
\]
for $\omega\in C^n(\calu_q,\calh),\ \ X_i\in \calu_q$. One has the 
following lemma.

\begin{lemma}
Let $\Psi\in \calh=C^0(\calu_q,\calh)$; the following
conditions $(i)$, $(ii)$ and $(iii)$ are equivalent\\
$(i)$ $d^k(\Psi)=0$ for some $k$ with $1\leq k\leq N-1$\\
$(ii)$ $\Psi\in\calh_I$\\
$(iii)$ $d^n(\Psi)=0$ for any $n\in\{1,\dots,N-1\}$.
\end{lemma}
Observe first that $d(=d_1)$ coincides in degree 0 with the 
Hochschild differential. Then the result is a consequence of the 
following formula which one proves by induction on $n$.
 
\[
d^n\Psi(\bbbone,\dots,\bbbone,X)=(1+q^2)\dots
(1+q^2+\dots+q^{2(n-1)})d\Psi(X)
\]
 for $\Psi\in C^0(\calu_q,\calh)$ and  for any $n\geq 1$, $X\in 
\calu_q$ where $\bbbone$ is the unit of
$\calu_q$.\\

This lemma implies : 
$H^0_{(k)}(C(\calu_q,\calh),d)=H^0(\calu_q,\calh)$, $\forall
k\in\{1,\dots,N-1\}$. This is a special case of Theorem 2 of Section 
6.
As an easy consequence, one obtains the following
result.
\begin{proposition}
The $N$-complex $(\calh^\bullet,d)$ can be canonically
identified with the $N$-subcomplex of $(C(\calu_q,\calh),d)$
generated by $\calh$. 
\end{proposition}

Thus one has $\calh^\bullet \subset C(\calu_q,\calh)$ and the
$N$-differential $d$ of $C(\calu_q,\calh)$ extends the one of
$\calh^\bullet$; we now extend $A$ to $C(\calu_q,\calh)$.
 
\begin{lemma}
Let us extend $A$ to $C(\calu_q,\calh)$ as a homogeneous
endomorphism $\omega\mapsto (A\omega)$ of degree $0$ by setting
\[
(A\omega)(X_1,\dots,X_n)=q^{2n}A\omega(X_1,\dots,X_n)
\]
for $\omega\in C^n(\calu_q,\calh)$ and $X_i\in\calu_q$. On 
$C(\calu_q,\calh)$one has
 $Ad-q^2dA=0$, $A^N=0$ and $(d+A)^N=0$.
 \end{lemma}
  
We have now extended to $C(\calu_q,\calh)$ the whole structure
defined previously on $\calh^\bullet$. Indeed the uniqueness in Lemma 
14 implies
that $A$ defined on $C(\calu_q,\calh)$ in last lemma is an
extension of $A$ defined on $\calh^\bullet$ in Lemma 14.
One then extends to $C(\calu_q,\calh)$ the definition of $Q$ by
setting again $Q=d+A$.\\

As explained in Section 6, Theorem 2 (1), the spaces
$H^n_{(k)}(C(\calu_q,\calh),d)$ can be computed in terms of the
Hochschild cohomology $H(\calu_q,\calh)$. In particular, one sees
that $H^n_{(k)}(C(\calu_q,\calh),d)$ does not generally vanish
for $n\geq 1$. This implies that one cannot expect for the
generalized homology of $Q$ on $C(\calu_q,\calh)$ such a simple
result as the one given by Theorem 5 for the generalized homology
of $Q$ on $\calh^\bullet$. Nevertheless, in view of Lemma 15, one
has
$H^0_{(k)}(C(\calu_q,\calh),d)=\calh_I=H^0_{(k)}(\calh^\bullet,d)$
and therefore one may expect 
$H^0_{(k)}(C(\calu_q,\calh),Q)=H_{(k)}(\calh_I,A)(=H_{(k)}
(\calh^\bullet,Q))$. In fact, this is essentially true. However
some care must be taken because $Q$ is not homogeneous so
$H_{(k)}(C(\calu_q,\calh),Q)$ is not a graded vector space.
Instead of a graduation, one has an increasing filtration
$F^nH_{(k)}(C(\calu_q,\calh),Q)$, ($n\in \mathbb Z$),
with
$F^nH_{(k)}(C(\calu_q,\calh),Q)=0$ for $n<0$ and where, for 
$n\geq
0$, $F^nH_{(k)}(C(\calu_q,\calh),Q)$ is the canonical image  in
$H_{(k)}(C(\calu_q,\calh),Q)$ of $\ker(Q^k)\cap
\displaystyle{\oplusinf^{r=n}_{r=0}}C^r(\calu_q,\calh)$. There is
an associated graded vector space
\[
^{\mathrm gr}H_{(k)}(C(\calu_q,\calh),Q)=\oplusinf_n
F^nH_{(k)}(C(\calu_q,\calh),Q)/F^{n-1}H_{(k)}(C(\calu_q,\calh),Q)
\]
which here is $\mathbb N$-graded. One has
$F^0H_{(k)}(C(\calu_q,\calh),Q)=^{\mathrm
gr}H_{(k)}^0(C(\calu_q,\calh),Q)$
and it is this space which is the correct version of the
$H^0_{(k)}(C(\calu_q,\calh),Q)$ above in order to identify
$H_{(k)}(\calh_I,A)$ in the generalized homology of $Q$ on
$C(\calu_q,\calh)$.
\begin{theo}
The inclusion $\calh^\bullet\subset C(\calu_q,\calh)$ induces the
 isomorphisms 
\[
H_{(k)}(\calh^\bullet,Q)\simeq
F^0H_{(k)}(C(\calu_q,\calh),Q)\  \mathrm{for}\   1\leq k\leq N-1.
\]
In
particular, with obvious identifications, one has
\[
F^0H_{(k)}(C(\calu_q,\calh),Q)=H_{(k)}(\calh_I,A),\ \ \ \forall
k\in \{1,\dots,N-1\}.
\]
\end{theo}
The proof is not difficult, for it as well as for complete proofs of 
all the results of this section we refer to \cite{D-VT2}.\\

If one compares this construction involving Hochschild cochains 
with the preceeding one, what has
been gained here besides the explicit occurrence of the quantum
gauge aspect is that the extended space $C(\calu_q,\calh)$ is a
tensor product $\calh\otimes\calh'$ of the original space
$\calh$ with the tensor algebra $\calh'=T(\calu^\ast_q)$ of the
dual space of $\calu_q$. The factor $\calh'$ can thus be
interpreted as the state space for some generalized ghost.
What has been lost is the minimality of the generalized
homology, i.e. besides the ``physical" $H_{(k)}(\calh_I,A)$, the
generalized homology of $Q$ on $\calh\otimes \calh'$ contains
some other non trivial subspace in contrast to what happens on
$\calh^\bullet$. In the usual homological (BRS) methods however
such a ``non minimality" also occurs. Indeed,as explained in last 
section, in
the homological approach to constrained classical systems, the
relevant homology contains besides the functions on the reduced
phase space the whole cohomology of longitudinal forms. The same is 
true for the BRS cohomology of gauge
theory \cite{BRS}, \cite{bc-r}.\\
In the usual situations where one applies the BRS construction
(gauge theory, constrained systems) one has a Lie algebra
$\fracg$ (the Lie algebra of infinitesimal gauge
transformations) acting on some space $\calh$ and what is really
relevant at this stage is the Lie algebra cohomology 
$H(\fracg,\calh)$ of
$\fracg$ acting on $\calh$. The extended space is then the space
of $\calh$-valued Lie algebra cochains of $\fracg$,
$C(\fracg,\calh)$. This extended space is thus also a tensor
product $\calh\otimes\calh'$ but now $\calh'$ is the exterior
algebra $\calh'=\Lambda\fracg^\ast$ of the dual space of
$\fracg$. That is why this factor can be
interpreted (due to antisymmetry) as a fermionic state space; 
indeed that is the reason why one
gives a fermionic character to the ghost \cite{BRS},
\cite{bc-r}. There is however another way
to proceed in these situations which is closer to what has been
done in our case here. To understand it, we recall that any
representation of $\fracg$ in $\calh$ is also a  
representation of the enveloping algebra $U(\fracg)$ in $\calh$.
Thus $\calh$ is a left $U(\fracg)$-module. Since $U(\fracg)$ is a
Hopf algebra, one can convert as above $\calh$ into a bimodule for
$U(\fracg)$ by taking as right action the trivial representation
given by the counit. It turns out that as explained in Section 3, 
Theorem 1, the $\calh$-valued Hochschild
cohomology of $U(\fracg)$, $H(U(\fracg),\calh)$, coincides with
the $\calh$-valued Lie algebra cohomology of $\fracg$,
$H(\fracg,\calh)$.
Since it is the latter space which is relevant one can as well
take as extended space the space of $\calh$-valued Hochschild
cochains of $U(\fracg)$, $C(U(\fracg),\calh)$, and then compute
its cohomology. Again this space is a tensor product
$\calh\otimes \calh'$ but now $\calh'=T(U(\fracg)^\ast)$ is a
tensor algebra as in our case.

\section*{Acknowledgments}

I thank Robert Coquereaux and Raymond Stora for their constructive 
critical reading of the manuscript.

\appendix

\section{Remarks on tensor products}

Let $\wedge\{d\}$ be the associative unital $\kb$-algebra generated 
by an element $d$ satisfying $d^2=0$. As a $\kb$-algebra 
$\wedge\{d\}=\kb\bbbone \oplus \kb d$ is the exterior algebra (over 
$\kb$) of the free $\kb$-module of rank one. It is clear that a 
$\wedge\{d\}$-module is the same thing as a differential module (as 
defined in Section 2). Given two differential modules $E$ and $F$ 
there is a canonical structure of $\wedge\{d\}\otimes 
\wedge\{d\}$-module on $E\otimes F$, where the first factor (resp. 
the second factor) corresponds to the structure  of 
$\wedge\{d\}$-module of $E$ (resp. of $F$). To say that for any such 
$E$ and $F$ there is a canonical differential on $E\otimes F$ (i.e. a 
canonical structure of $\wedge\{d\}$-module on $E\otimes F$) which 
only depends on the differentials of $E$ and $F$ (i.e. on their 
$\wedge\{d\}$-module structures) is the same thing as to say that one 
has a coproduct $\Delta$ on $\wedge\{d\}$, that is a homomorphism of 
unital $\kb$-algebras $\Delta:\wedge\{d\}\rightarrow 
\wedge\{d\}\otimes \wedge \{d\}$. One must have 
$\Delta(\bbbone)=\bbbone\otimes\bbbone$ so $\Delta$ is fixed by 
giving a $\Delta(d)\in\wedge\{d\}\otimes \wedge \{d\}$ satisfying 
$(\Delta(d))^2=0$. One has $\Delta(d)=\alpha\bbbone\otimes 
\bbbone+\beta\bbbone\otimes d+\gamma d\otimes \bbbone +\delta 
d\otimes d$ with $\alpha, \beta, \gamma, \delta\in \kb$ and 
$(\Delta(d))^2=0$ implies $\alpha^2=0,2\alpha\beta=0,2\alpha\gamma=0$ 
and $2(\alpha\delta+\beta\gamma)=0$.\\

Let us now assume that $\kb$ is a field of characteristic different 
from 2. Then the above conditions imply $\alpha=0$ and $\beta\gamma=0$ 
i.e. either $\Delta(d)=\beta\bbbone \otimes d+\delta d\otimes d$ or 
$\Delta(d)=\gamma d\otimes \bbbone+\delta d\otimes d$. It is already 
clear that generically the differential $\beta\bbbone\otimes d+\delta 
d\otimes d$ (resp. $\gamma d\otimes \bbbone+\delta d\otimes d$) on 
$E\otimes F$ will lead to a homology $H(E\otimes F)$ for $E\otimes F$ 
different from $H(E)\otimes H(F)$. Notice that if one imposes the 
natural requirement of coassociativity for $\Delta$ one is led to the 
only 3 possibilities $\bbbone\otimes d$, $d\otimes \bbbone$ or 
$d\otimes d$ for the differential on the tensor products.\\

Let us come back to a general ring $\kb$. Consider the associative 
unital $\kb$-algebra $\cald_{-1}$ generated by two elements $d$ and 
$\Gamma$ satisfying $d^2=0$, $\Gamma^2=\bbbone$ and $\Gamma 
d=-d\Gamma$. This algebra is a Hopf algebra for the counit 
$\varepsilon$, the antipode $S$ and the coproduct $\Delta$  given by 
: $\varepsilon(d)=0$, $\varepsilon(\Gamma)=1$, $S(d)=-\Gamma d$, 
$S(\Gamma)=\Gamma$, $\Delta(d)=d\otimes \bbbone+\Gamma\otimes d$ and 
$\Delta(\Gamma)=\Gamma\otimes\Gamma$. The Hopf algebra $\cald_{-1}$ 
can be understood as a version of the universal enveloping algebra of 
the super Lie algebra with only one odd element $d$ such that 
$[d,d]=0$. Let $E=E^0\oplus E^1$ be a $\mathbb Z_2$-complex then $E$ 
is a $\cald_{-1}$-module if $d$ is represented by the differential of 
$E$ and if $\Gamma$ is represented by the multiplication by $(-1)^i$ 
on $E^i$ for $i\in \{0,1\}$. One verifies easily that the tensor 
product of $\mathbb Z_2$-complexes defined in Section 3 corresponds 
to the above structure, i.e. that it is induced by the coproduct 
$\Delta$. Thus one can understand the tensor product of complexes in 
terms of a Hopf algebra. We now show that the same is true for 
$N$-complexes.\\

Let $q\in \kb$ be such that Assumption ($A_1$) of Section 5 is 
satisfied and let us consider the associative unital $\kb$-algebra 
$\cald_q$ generated by two elements $d$ and $\Gamma$ satisfying 
$d^N=0$, $\Gamma^N=\bbbone$ and $\Gamma d=qd\Gamma $. Again $\cald_q$ 
is a Hopf algebra for the counit $\varepsilon$, the antipode $S$ and 
the coproduct $\Delta$   given by : $\varepsilon(d)=0$, 
$\varepsilon(\Gamma)=1$, $S(d)=-\Gamma^{N-1}d$, 
$S(\Gamma)=\Gamma^{N-1}$, $\Delta(d)=d\otimes \bbbone +\Gamma\otimes 
d$ and $\Delta \Gamma=\Gamma\otimes \Gamma$. Let 
$E=E^0\oplus\cdots\oplus E^{N-1}$ be a $\mathbb Z_N$-complex (see in 
Section 6) then $E$ is a $\cald_q$-module if $d$ is represented by 
the $N$-differential of $E$ and $\Gamma$ is represented by the 
multiplication by $q^i$ on $E^i$ for $i\in\{0,\dots,N-1\}$. Again one 
verifies easily that the $(q)$ tensor product of $\mathbb 
Z_N$-complexes (or of $N$-complexes) defined in Section 8 (introduced 
originally in \cite{Kap}) is induced by the coproduct of $\cald_q$.\\

Finally it is worth noticing here that instead of $\cald_{-1}$ one 
can use the exterior algebra of the free module of rank one 
$\wedge^\bullet\{d\}=\wedge\kb=\kb\bbbone\oplus\kb d$ considered as a 
graded Hopf algebra. That is, as an associative algebra 
$\wedge^\bullet\{d\}$ is isomorphic to $\wedge\{d\}$ but it is a 
$\mathbb Z_2$-graded algebra with $\wedge^0\{d\}=\kb\bbbone$, 
$\wedge^1\{d\}=\kb d$ and furthermore it is a graded Hopf algebra for 
the counit $\varepsilon$, the antipode $S$ and the coproduct $\Delta$ 
given by $\varepsilon(d)=0$, $S(d)=-d$ and $\Delta 
d=d\otimes\bbbone+\bbbone\otimes d$ where now $\Delta$ is a 
homomorphism of graded algebras of $\wedge^\bullet\{d\}$ into 
$\wedge^\bullet\{d\}\otimes \wedge^\bullet\{d\}$ with 
$\wedge^\bullet\{d\}\otimes \wedge^\bullet\{d\}$ being the (twisted) 
tensor product of graded algebras defined in Section 8. A $\mathbb 
Z_2$-complex is canonically a graded $\wedge^\bullet\{d\}$-module and 
the tensor product of complexes can be also defined by using the 
above graded coproduct. Of course $\wedge^\bullet\{d\}$ is also a 
version (which is graded) of the universal enveloping algebra of the 
super Lie algebra with one odd generator $d$ satisfying $[d,d]=0$. 
The advantage of $\cald_{-1}$ is that it generalizes as $\cald_q$ for 
$N$-complexes as explained above and that it is an ordinary Hopf 
algebra (in fact  $\Gamma$ plays the role of the graduation).

\section{Longitudinal forms}

Let $\cala$ be an associative unital algebra over $\mathbb R$ or 
$\mathbb C$ (here $\kb=\mathbb R$ or $\mathbb C$) and let us denote 
by $Z(\cala)$ the center of $\cala$ that is the commutative unital 
subalgebra of $\cala$ defined by $Z(\cala)=\{z\in\cala\vert za=az,\ \ 
\forall a\in \cala\}$. One has $\cala=Z(\cala)$ if and only if 
$\cala$ is commutative. Recall that a derivation of $\cala$ is a 
linear mapping $X:\cala\rightarrow \cala$ such that one has (Leibniz 
rule) $X(ab)=X(a)b+aX(b)$, $\forall a,b\in \cala$. If $X$ and $Y$ are 
derivations of $\cala$, their composition $XY$ (product in 
$\mathrm{End} (\cala)$) is not a derivation but the commutator 
$[X,Y]=XY-YX$ is again a derivation of $\cala$.  On the other hand, 
if $X$ is a derivation of $\cala$ and if $z$ is in the center of 
$\cala$ then $zX$ (defined by $(zX)(a)=zX(a),\ \ \forall a\in\cala$) 
is again a derivation of $\cala$. Thus the vector space $\Der(\cala)$ 
of all derivations of $\cala$ is a Lie algebra (for $[\cdot,\cdot]$) 
and also a $Z(\cala)$-module, both structures being connected through 
the identity\\
$(DZ) \hspace{2cm}  [X,zY]=z[X,Y]+X(z)Y$\\
for any derivation $X$ and $Y$ and for any $z\in Z(\cala)$; one 
verifies easily that the center is stable by derivation i.e. that one 
has $X(z)\in Z(\cala)$ for any $X\in\Der(\cala)$ and $z\in Z(\cala)$. 
Let $\calf\subset\Der(\cala)$ be a $Z(\cala)$-submodule which is also 
a Lie subalgebra of $\Der(\cala)$. The graded space 
$C_\wedge(\calf,\cala)$ of (Chevalley-Eilenberg) $\cala$-valued 
cochains of the Lie algebra $\calf$ (see in Section 3) is canonically 
a graded algebra and, since $\calf$ operates by derivation on 
$\cala$, the corresponding Chevalley-Eilenberg differential $d$ is an 
antiderivation of $C_\wedge(\calf,\cala)$. Thus 
$C_\wedge(\calf,\cala)$ is a graded differential algebra (see in 
Section 8). It follows from the above identity $(DZ)$ that the graded 
subalgebra $\os_\calf(\cala)$ of cochains which are 
$Z(\cala)$-multilinear is stable by $d$ so $\os_\calf(\cala)$ is a 
graded  differential algebra.\\

Let $V$ be a smooth manifold and let $\calf$ be a foliation of $V$ 
(see in Section~9), then the graded differential algebra 
$\os_\calf(C^\infty(V))$ is referred to as the graded differential 
algebra of {\sl longitudinal forms} and is denoted by 
$\Omega(V,\calf)$; its elements are called {\sl longitudinal forms}. 
This is the graded differential algebra considered in Section 9. 
Notice that when $\calf$ coincides with the module 
$\Der(C^\infty(V))$ of all vector fields on $V$ then 
$\os_\der(C^\infty(V))=\os_{\der(C^\infty(V))}(C^\infty(V))$ is the 
graded differential algebra $\Omega(V)$ of differential forms on $V$. 
This is why $\os_{\derth}(\cala)=\os_{\derth(\cala)}(\cala)$ is a 
noncommutative generalization of the graded differential algebra of 
differential forms when $\cala$ is noncommutative; there are other 
noncommutative generalizations of differential forms (see e.g. in 
\cite{D-V4}).

\end{document}